\documentclass[reqno]{amsart}
\usepackage{bbm}
\usepackage{mathrsfs}
\usepackage{amsmath}
\allowdisplaybreaks[0]
\usepackage{paralist}
\usepackage{graphics}
\usepackage{epsfig}
\usepackage[pdftex]{hyperref}
\usepackage{amsfonts}
\usepackage{CJK}
\usepackage{fancyhdr}
\usepackage{graphicx}
\usepackage[dvipsnames,usenames]{color}
\usepackage{titletoc}
\usepackage{latexsym}
\usepackage{amssymb}
\usepackage{multicol}
\usepackage{graphics}
\usepackage{subfigure}
\usepackage{indentfirst}
\usepackage{cases}
\usepackage{curves}
\usepackage{cite}
\newtheorem{theorem}{Theorem}[section]

\newtheorem{lemma}{Lemma}[section]
\newtheorem{definition}{Definition}[section]

\newtheorem{remark}{Remark}[section]
\newtheorem{example}{Example}[section]

\numberwithin{equation}{section}

\def\x#1{(\ref{#1})}
\def\R{{\Bbb R}}

\def\N{{\Bbb N}}

\def\E{{\mathbb{E}}}

\def\bc{\begin{center}}
\def\ec{\end{center}}
\def\ba{\begin{array}}
\def\ea{\end{array}}
\def\be{\begin{equation}}
\def\ee{\end{equation}}
\def\bea{\begin{eqnarray}}
\def\eea{\end{eqnarray}}
\def\beaa{\begin{eqnarray*}}
\def\eeaa{\end{eqnarray*}}

\def\ben{\begin{enumerate}}
\def\een{\end{enumerate}}
\def\hh{\!\!\!\!}
\def\EM{\hh &   &\hh}
\def\EQ{\hh & = & \hh}

\def\LE{\hh & \le & \hh}
\def\GE{\hh & \ge & \hh}

\def\nn{\nonumber}
\def\oo{\infty}
\def\ifl{\iffalse}
\def\lb{\label}
\def\prf{\mbox{\bf Proof.~}}

\DeclareMathOperator{\im}{im}

\title[]{Nonuniform Mean-square Exponential Dichotomies and Mean-square Exponential  Stability}

\author[Hailong Zhu, Li Chen, Xiuli He]
{Hailong Zhu $^{1,2}$, \quad  Li Chen $^{2}$, \quad Xiuli He$^{3}$}
\address{$^1$ Anhui
University of Finance and Economics, Bengbu, 233030, China}
\address{$^2$ Universit\"{a}t Mannheim, 68131, Mannheim, Germany}
\address{$^3$  Department of Mahthematics,  Hohai University, Nanjing, 210098, China}
\email{hai-long-zhu@163.com (H. Zhu)}
\email{chen@math.uni-mannheim.de (L. Chen)}
\email{hexiu00@163.com (X. He)}

\begin{document}
\setlength{\parskip}{0.5\baselineskip}

\begin{abstract} In this paper, the existence conditions of nonuniform  mean-square exponential   dichotomy (NMS-ED)  for a linear stochastic differential equation (SDE) are established.
The difference of the conditions for the existence of a nonuniform dichotomy between an SDE and an ordinary differential equation (ODE) is that the first one needs an additional assumption, nonuniform Lyapunov matrix, to guarantee that the linear
SDE  can be transformed into a decoupled one,  while the second does not.
Therefore,
the first main novelty of our work is that we establish some preliminary results  to tackle  the stochasticity.
This paper is also  concerned with  the  mean-square exponential stability of nonlinear perturbation of a linear SDE under the condition of nonuniform  mean-square exponential contraction  (NMS-EC).
For this purpose,   the concept of  second-moment regularity coefficient is introduced. This concept
 is essential in determining  the stability
of the perturbed equation,
and hence we deduce the lower and upper  bounds of this coefficient.
 Our results imply that the lower and upper  bounds of the second-moment regularity
 coefficient  can be expressed solely by the drift term of the linear SDE.
\end{abstract}

\maketitle

\section{\bf Introduction}
\setcounter{equation}{0} \noindent

Mean-square dynamical behavior is one of  the important
concepts to describe the flows produced by SDEs  or random differential equations (RDEs). This is due to the fact that in the case of mean-square setting, the dynamical behavior of SDEs and RDEs  are essentially
deterministic with the stochasticity built into or hidden in the
time-dependent state spaces (under specific conditions,  there is no difference between the flows generated by the SDEs and RDEs;  in fact, the flow of SDEs is conjugate to the flow of RDEs \cite{il-01}). Over   the   years,   its  many  properties  and   corresponding results have been presented by many researchers.
For example, Kloeden and Lorenz \cite{kl}
provided a definition of mean-square random dynamical systems
and studied the existence of pullback attractors.
%(see e.g., \cite{ar98} for details about random dynamical systems).
In \cite{fl-10, ls, zcz}, the concept of mean-square almost
automorphy for stochastic process was introduced,   the existence, uniqueness and asymptotic stability
of mean-square almost automorphic solutions of SDEs were established respectively.   Using a stochastic version of theta
method, Higham \cite{hig} combined analytical and numerical techniques to tackle mean-square asymptotic stability for SDEs. Recently,  Zhu and Chu \cite{zc1}
 presented the numerical methods for
 a mean-square exponential dichotomy (MS-ED)  of a linear SDE and showed that  the MS-ED
is equivalent to  the numerical results for sufficient small step
sizes under natural conditions. We also refer to  \cite{drk, hms, hmy, zc2} for more related results and techniques about this topic.

The concept of MS-ED  is extended from the classical notation of exponential dichotomy, which can be traced back to Perron \cite{per} in 1930s. Since then it has become a very important part of the general theory of
dynamical systems, particularly in what concerns the study of stable and
unstable invariant manifolds, and therefore has attracted much attention during
the last few decades.  One can see, for example, \cite{huy-06,  pa-88, pp-06, rr-95, lrs, ss} about evolution equations, \cite{cof-71,lin, per-71} about functional differential   equations,
  \cite{cl-94,cl,lmr,ss1} about skew-product   flows, and  \cite{drk,st, zlz, zc2, zcz} about random  systems  or   stochastic  equations. We also refer to the books \cite{chi, cop, ms} for details and further references
related to exponential dichotomies.

%Stochastic differential equations (SDEs) have been studied from different angles, since SDEs have important applications in various scientific areas. For example, Arnold \cite{ar74} and Mao \cite{mao} study and discuss the theory and applications of SDEs. Kloeden and Platen \cite{kp}  summarize the development of the theory of the numerical integration of SDEs.  Pardoux and R\v{a}\c{s}canu \cite{pr} consider the backward stochastic differential equations.
%For more topics about  SDEs , one can see \cite{ev, ga, lad, O-13} and
%the   references  therein.
%
%Among those topics,

However, dynamical  systems  exhibit  various  different   kinds  of   dichotomic    behaviors, and the notion of classical exponential dichotomy cannot contain all possible dichotomic behaviors, as   Barreira  and   Valls  mentioned  in \cite{bv06}, ``the notion of exponential dichotomy demands considerably from
the dynamics and it is of considerable interest to look for more general types of
hyperbolic behavior".   In these years, many attempts have been made (see, e.g.,\cite{np,np-97, pp-04}) to extend the concept of the classical  dichotomies. For more recent works we mention in particular the papers \cite{bcv11,bcv13,bv05,bv06,bv07,bv08,bv10}, which,
inspired by the fundamental work of nonuniformly  hyperbolic trajectory introduced in \cite{bp02,bp07}, extend  the  concept  of   exponential   dichotomy to the nonuniform  ones and investigate some related problems.
In fact, exponential  dichotomy implies nonuniform exponential  dichotomy (see e.g., \cite{bv06,bv07,bv08}). However, the contrary is not true in general.
For example, Barreira  and   Valls \cite{bv07} showed that the linear equation
\[u'=(-a-bt \sin t)u,\quad v'=(a+bt \sin t)v\]
with $a>b>0$ admits   a   nonuniform  exponential  dichotomy  but  does   not  admit   a   uniform  exponential
dichotomy.

As our knowledge, the concept of
MS-ED was first introduced by Stanzhyts'kyi \cite{sk-01}, in which a sufficient condition has been proved to ensure that a linear SDE admits an MS-ED. Based on the definition of MS-ED, Stanzhyts'kyi and Krenevych \cite{sk-06}  proved the existence of a
quadratic form of the linear SDE. In \cite{zc2} the robustness of MS-ED for a linear SDE was established, and Stoica \cite{st} studied  stochastic
cocycles  in Hilbert spaces by using MS-ED. Recently, Doan et al. \cite{drk} considered the MS-ED spectrum for random dynamical system.

 Now we recall the definition of MS-ED.
 Consider the following linear $n$-dimensional It\^{o} stochastic system
\be\lb{a1}du(t)=A(t)u(t)dt+G(t)u(t)d\omega(t), \quad t\in I,\ee
where $I$  is  either  the  half line $\R^{+}$ or the  whole line $\R$, and $A(t)=(A_{ij}(t))_{n \times n}$, $G(t)=(G_{ij}(t))_{n \times n}$
are continuous functions with real entries, which satisfy
 \be\lb{a4} \limsup_{t\rightarrow +\oo}\log^{+}\|A(t)\|=0,\quad{\rm and}\quad \limsup_{t\rightarrow +\oo}\log^{+}\|G(t)\|=0,
\ee
with $\log^{+}x=\max\{0, \log x\}$.
 Eq. \x{a1} is said to possess a \emph{mean-square exponential dichotomy } if there exist linear projections $P(t): L^{2}(\Omega, \R^{n}) \rightarrow L^{2}(\Omega, \R^{n})$ such that
\be\lb{c1}\Phi(t)\Phi^{-1}(s)P(s)=P(t)\Phi(t)\Phi^{-1}(s),  \quad
\forall~ t,s\in I,\ee and positive constants $K, \alpha$  such that  \be\lb{a5}
\begin{split}
\E\|\Phi(t)\Phi^{-1}(s)P(s)\|^2 \leq K e^{-\alpha(t-s)},
\quad \forall~ (t,s)\in I^{2}_{\geq},\\
\E\|\Phi(t)\Phi^{-1}(s)Q(s)\|^2 \leq K e^{-\alpha(s-t)}, \quad
\forall~ (t,s)\in I^{2}_{\leq},
\end{split}
\ee
where  $\Phi(t)$ is a fundamental matrix solution of {\rm\x{a1}}, and $Q(t)={\rm Id}-P(t)$ is the
complementary projection of $P(t)$ for each $t\in I$.
$I^{2}_{\geq}:=\{(t,s)\in I^{2}: t\geq s\}$ and
$I^{2}_{\leq}:=\{(t,s)\in I^{2}: t\leq s\}$ denote the relations of $s$ and $t$ on $I$. The  constants $\alpha$ and  $K$ are called
the   \emph{exponent} and   the   \emph{bound} respectively  in  the case of  deterministic systems \cite{hen-81}.

This paper, inspired by both the mean-square dynamical properties and the nonuniform behavior, is to study the NMS-ED and its related problems. Eq.  \x{a1} is said to possess a  \emph{nonuniform mean-square exponential dichotomy} if  there exist linear projections $P(t): L^{2}(\Omega, \R^{n}) \rightarrow L^{2}(\Omega, \R^{n})$  such that \x{c1} holds,
and positive constants $K, \alpha$ and $\beta\in [0, \alpha)$ such that  \be\lb{c2}
\begin{split}
\E\|\Phi(t)\Phi^{-1}(s)P(s)\|^2 \leq K e^{-\alpha(t-s)+\beta s},
\quad \forall~ (t,s)\in I^{2}_{\geq},\\
\E\|\Phi(t)\Phi^{-1}(s)Q(s)\|^2 \leq K e^{-\alpha(s-t)+\beta s}, \quad
\forall~ (t,s)\in I^{2}_{\leq}.
\end{split}
\ee
 For   the convenience of statement, in the rest of this paper,  we  call  $\alpha$ the \emph{exponent}, $K$ the \emph{bound}, and $\varepsilon$ the \emph{nonuniform degree}. From the point of dichotomic    behavior, the standard growth condition \x{a5} on $\Phi$ is replaced by a
much weaker condition \x{c2} so that the main results can be applied to
a larger class of equations.
The  nonuniformity in \x{c2}   indicates that the bound of  the corresponding solution
depends   on  initial time $s$ (while in the uniform case \x{a5} this bound must be chosen
independently of $s$).
Clearly,  if one considers $\beta=0$ in \x{c2}, we say that \x{a1} admits a (uniform) mean-square exponential dichotomy \x{a5}. That is to say,  a mean-square exponential dichotomy is a particular  case  of the nonuniform ones.
On the contrary, the nonuniform part $e^{\beta s}$ in \x{c2} cannot be removed in some cases.  For example, let $a>b>0$ be real parameters,
  \[
\left\{ \begin{array}{ll}
du & =(-a-bt\sin t)u(t)dt+\sqrt{2b\cos t}\exp(-at+bt\cos t)d\omega(t)\\
dv & =(a+bt\sin t)v(t)dt-\sqrt{2b\cos t}\exp(at-bt\cos t)d\omega(t)
\end{array} \right.
\]
admits an NMS-ED  which is not uniform. See Example 6.1 in \cite{zz} for details.

The first aim of this paper is: under which conditions the
NMS-ED  of \x{a1} exists? In the process of establishing the existence conditions of nonuniformity, a significant difference between ODEs  and SDEs
can be observed,
%For SDE \x{a1}, the main difficulty is caused by diffusion term $G(t)$, which should be %tackled by other techniques as below.
that is, for an ODE $x'=A(t)x$, one can assume that $A(t)$ has the block form
 \[A(t)=\left(\begin{array}{ll}
A_{1}(t) & 0\\
0 & A_{2}(t)
\end{array}\right).\]
The blocks $A_{1}(t)$, $A_{2}(t)$ correspond, respectively, to stable and unstable components of $A(t)$, under which the system $x'=A(t)x$ can be proved to have a nonuniform exponential dichotomy (see \cite{bv07} for details).
However, this assumption cannot be used  directly for SDE \x{a1}, since it is   unreasonable to assume that $A(t)$ and $G(t)$ in system \x{a1} can be decoupled into block forms with the same dimensions.
To overcome the difficulty caused by the fact that block forms $A(t)$ and $G(t)$ may have different dimensions, a condition called nonuniform Lyapunov matrix is introduced, under which \x{a1} can be transformed into a new system
%\[dv(t)=B(t)v(t)dt+v(t)d\omega(t).\]
\bea\lb{c15}dv(t)\EQ B(t)v(t)dt+v(t)d\omega(t)\nn\\
\EQ\left(\begin{array}{ll}
B_{1}(t) & 0\\
0 & B_{2}(t)
\end{array}\right)v(t)dt+\left(\begin{array}{ll}
Id_{n_{1}\times n_{1}} & 0_{n_{1}\times n_{2}}\\
0_{n_{2}\times n_{1}} & Id_{n_{2}\times n_{2}}
\end{array}\right)v(t)d\omega(t).\eea
Thus the drawback in stochasticity  can be overcome since the unit matrix can be seen as a block form.
 %Thus the assertion for \x{a1} admits a NMS-ED  can be proved by using nonuniform Lyapunov matrix $S(t)$.

\begin{theorem}\lb{main1}Assume that  there is a nonuniform Lyapunov matrix $S(t)$, which transforms {\rm \x{a1}} into the block form {\rm \x{c15}}.
 Then for sufficiently small $\varepsilon>0$,  {\rm \x{a1}} admits an
 NMS-ED  with  the exponent
\[\alpha =\max\{-(\chi_{k}+\varepsilon), \chi_{k+1}+\varepsilon\}>0,\]
where the notations $\chi_{k}, k=1,\ldots,r (r\le n)$ are the second-moment Lyapunov exponents given in {\rm\x{b10}}.
\end{theorem}

 Theorem \ref{main1} is   on  the existence of
 NMS-ED of system \x{a1}, which   is   a  generalization  of  nonuniform dichotomy for ODEs. The proof of Theorem \ref{main1} is presented in Section 3, which is much more delicate than that of  previous works for ODEs (see \cite{bv07}). In fact, a linear SDE which is nonuniformly kinematically similar to \x{a1} is constructed by nonuniform Lyapunov matrix,
 whereby several results are needed before the proof of Theorem \ref{main1}.

Next  we consider a nonlinear SDE
\be\lb{a3} du(t)=\big{(}A(t)u(t)+f(t,u(t))\big{)}dt+\big{(}G(t)u(t)+h(t,u(t))\big{)}d\omega(t), \quad t\in I,
\ee
which is a perturbation of \x{a1}.
%we present a sufficient condition for the perturbation of mean-square
%exponentially stability.
The trivial solution of \x{a1} is said to be
\emph{mean-square exponentially stable} (or \emph{second-moment exponentially stable}) if there exist positive constants $C$, $\chi$ such that
\beaa\E\|x(t)\|^{2}\le C\|x_{0}\|^{2}e^{-\chi(t-t_{0})},\quad \forall t\ge t_{0}\eeaa
for all $x_{0}\in \R^{n}$.
It is well-known that mean-square exponential stability  is a
special case of $p$th moment exponential stability. This stability is
one of the most effective tools (for example, stability in probability,
moment stability and almost sure stability) to describe the stochastic
stability (see, e.g., \cite{ar74,mao-94,mao} for details), and
 mean-square exponential stability for  SDE  can be seen as a natural
 generalization of the classical concept of  exponential stability for ODEs
 (see e.g., \cite{cop}) since the It\^{o} stochastic calculus is a mean-square calculus.

 The second aim of this paper is to study the mean-square exponential stability of \x{a3} when \x{a1} admits
an NMS-EC, which is a special case of
NMS-ED with $P(t)=Id$ (see Section 4 for details).
Roughly speaking, NMS-EC of \x{a1} determines whether or not the trivial solution of the perturbed equation \x{a3} is mean-square exponential stability. Example \ref{exam61} in Section 6 indicates that in general the answer is negative. For ODEs, Lyapunov introduced regularity conditions to guarantee exponential stability of the trivial solution of the corresponding perturbed equation (see, e.g., \cite{bp02,bvgn}).  In order to generalize the Lyapunov stability theorem on the well-established deterministic theory, the notion of  \emph{regular} is  stated in the next section. Based on this additional assumption, NMS-EC indeed implies the stability of the trivial solution of \x{a3}.

\begin{theorem}\lb{main3}
  Assume that Eq. {\rm\x{a1}} admits a nonuniform mean-square exponential contraction, with the second-moment Lyapunov exponent $\chi$ of Eq. {\rm\x{a1}} being  regular. Then the trivial solution of the perturbed equation {\rm\x{a3}} is mean-square exponentially stable.
\end{theorem}

In addition, we draw this conclusion with a weaker hypothesis  in the following theorem. Roughly speaking, we obtain the mean-square exponential stability of the perturbed equation \x{a3}, which does not need the condition that Eq. \x{a1} is regular.

\begin{theorem}\lb{main2}Assume that  Eq. {\rm\x{a1}} admits a nonuniform mean-square exponential contraction with
$-q\alpha+\beta<0$ (see {\rm\x{e14}} and {\rm\x{a8}} in Section 4 for notations and details).
Then there exists $\delta>0$ sufficiently small so that for every initial condition $\xi_{0}\in \R^{n}$ with $\|\xi_{0}\|\le \delta$,
the solution of Eq. {\rm\x{a3}} starting at $\xi_{0}$  is   mean-square
exponential stable which satisfies:
\be\lb{e11}\E\|u(t)\|^{2}\le \tilde{K}e^{-\alpha t},\ee
where $\tilde{K}>0$ is a constant.
\end{theorem}
In Section 4, we start by proving this weaker statement.
After the proof of Theorem \ref{main2}, the fact that Theorem \ref{main3} can be obtained directly
from Theorem \ref{main2} is explained  in Remark
\ref{rem41}. In addition,
one can find that the second-moment regularity coefficient
$\gamma(\chi,\tilde{\chi})$ plays a key role in determining the stability
of the perturbed equation \x{a3} from the discussion of Remark \ref{rem41}. Hence,
our aim is to derive the lower and upper bounds of $\gamma(\chi,\tilde{\chi})$ in Section 5.

The paper is organized as follows.  The
next section introduces some notations and prepares several preliminary
results which will be used in later sections. Section 3 proves that \x{a1}
admits an NMS-ED by using nonuniform Lyapunov matrix $S(t)$.   Section 4 devotes to
the study of the mean-square exponential stability of  \x{a1}. Section 5
investigates the lower and upper bounds of the second-moment regularity coefficient
$\gamma(\chi,\tilde{\chi})$.  Finally, an example  is given in Section 6,
which shows that in general NMS-EC
is not enough to guarantee the stability of the perturbed equation of a linear SDE.
%The second one
%indicates that there esists a linear SDE which admits an NMS-ED but not uniform.

\section{Second-moment Lyapunov exponent}

\setcounter{equation}{0} \noindent

Thoughout this paper, we assume that $(\Omega, \mathscr{F},
\mathbb{P})$ is a  probability space,
$\omega(t)=(\omega_{1}(t),\ldots \omega_{n}(t))^{T}$ is an
$n$-dimensional Brownian motion defined on the space $(\Omega,
\mathscr{F}, \mathbb{P})$. $\|\cdot\|$ is used to stand for either  the Euclidean vector norm or the matrix norm  as appropriate, and $L^{2}(\Omega, \R^{n})$
represents  the space of all $\R^{n}$-valued random variables
$x: \Omega \rightarrow \R^{n}$ such that
$$\E \|x\|^2=\int_{\Omega}\|x\|^2 d \mathbb{P}<\oo.$$
For $x\in L^{2}(\Omega, \R^{n})$, let
\[\|x\|_{2}=\left(\int_{\Omega}\|x\|^2 d \mathbb{P}\right)^{1/2}.\]
Obviously, $L^{2}(\Omega, \R^{n})$ is a Banach space with the norm $\|x\|_{2}$.

Define the \emph{second-moment Lyapunov exponent} $\chi:\mathbb{R}^{n}\rightarrow \mathbb{R}\cup\{-\infty\}$ for a stochastic process $u: \R \rightarrow L^{2}(\Omega, \R^{n})$ by the formula
%The equation \x{dual1} is said to be \emph{kinematically similar} to the equation
%\be\lb{equa1}dv(t)=B(t)v(t)dt+v(t)d\omega(t)\ee of the same form if there exists a continuously differentiable invertible matrix $S(t)$ which satisfies the differential equation
%\[dS(t)=(A(t)S-SB(t)-S+G(t)S(t))dt+(G(t)S-S)d\omega(t).\]
%The change of variables $x=S(t)y$ then transforms \x{dual1} into \x{equa1}
\be\lb{b1}
\chi(u_{0})=\limsup_{t\rightarrow +\infty}\frac{1}{t}\log\E\Vert u(t)\Vert^{2},
\ee
where $u(t)$ is the solution of \x{a1} with the initial point $u(0)=u_{0}$. The uniqueness of the solution of \x{a1} for any given initial value is nicely described in the book by Mao \cite[Theorem 2.1, p. 93]{mao}. Thus it follows from the abstract theory of Lyapunov exponents (see e.g., \cite{bp02} for a detailed exposition)  that the function $\chi$ takes at most $r\leq n$ distinct values on $\mathbb{R}^{n}\backslash \{0\}$, say
\be \lb{b10}-\infty\leq\chi_{1}<\cdots<\chi_{k}<0\leq\chi_{k+1}<\cdots<\chi_{r}.\ee

Let $\Phi(t)$ be a fundamental matrix solution of \x{a1}. By \cite[Theorem 3.2.4]{mao},
$\Phi(t)$ is invertible with probability $1$ in $I$. To introduce the notion of  regularity for SDEs, we need the following lemma, which illustrates that the existence of the fundamental matrix solution of the adjoint equation of \x{a1}.

\begin{lemma}  \label{lem34} {\rm (see \cite[Theorem 2.3.1]{lad})} Let $\Phi(t)$ be a fundamental matrix solution of
{\rm\x{a1}}. Then $\Phi^{-1}(t)$ is a fundamental matrix
solution of the following stochastic differential equation
\be\lb{b5}d\tilde{u}(t)=\tilde{u}(t)[-A(t)+G^{2}(t)]dt-\tilde{u}(t)G(t)d\omega(t),\quad
t\in I.\ee
\end{lemma}

In fact,  Lemma \ref{lem34} can be verified by using It\^{o} product rule:
\[d(\Phi\Phi^{-1})=d\Phi\Phi^{-1}+\Phi d\Phi^{-1}+d\Phi d\Phi^{-1}=d 1 =0.\]
Clearly, $\Phi^{-T}(t)$ is a fundamental matrix
solution of the following SDE
\be\lb{b2}d\tilde{u}(t)=\left(-A(t)+G^{2}(t)\right)^{T}\tilde{u}(t)dt-G^{T}(t)\tilde{u}(t)d\omega(t)\quad
\ee due to \x{b5}, where $\left(-A(t)+G^{2}(t)\right)^{T}$ and $G^{T}(t)$ denote the transpose of $-A(t)+G^{2}(t)$ and $G(t)$ respectively. For \x{b2}, consider the associated second-moment Lyapunov exponent $\tilde{\chi}:\mathbb{R}^{n}\rightarrow \mathbb{R}\cup\{-\infty\}$ defined by
\be\lb{b3}
\tilde{\chi}(\tilde{u}_{0})=\limsup_{t\rightarrow +\infty}\frac{1}{t}\log\E\Vert \tilde{u}(t)\Vert^{2},
\ee
where $\tilde{u}(t)$ is the solution of \x{b2} with the initial value $\tilde{u}(0)=\tilde{u}_{0}$.  Again it follows from the abstract theory of Lyapunov exponents that $\tilde{\chi}$ can take at most $s\leq n$ distinct values on $\mathbb{R}^{n}\backslash \{0\}$, say $-\infty\leq\tilde{\chi}_{s}<\cdots<\tilde{\chi}_{1}$.

Now define the \emph{second-moment regularity coefficient} of $\chi$ and $\tilde{\chi}$ by
\be\lb{b4}\gamma(\chi,\tilde{\chi})=\min\max\{\chi(u_{i})+\tilde{\chi}(\tilde{u}_{i}):1\le i\le n\},\ee
where the minimum is taken over all bases $u_{1},\ldots, u_{n}$ and $\tilde{u}_{1},\ldots, \tilde{u}_{n}$
of $\R^{n}$ such that $\langle u_{i}, \tilde{u}_{j}\rangle=\delta_{ij}$ for each $i$ and $j$ (here $\delta_{ij}$ is the Kronecker symbol). We say that a basis $(u_{1},\ldots, u_{n})$ is \emph{dual} to  a basis $(\tilde{u}_{1},\ldots, \tilde{u}_{n})$ if $\langle u_{i}, \tilde{u}_{j}\rangle=\delta_{ij}$ for each $i$ and $j$. The second-moment Lyapunov  exponents $\chi$ and $\tilde{\chi}$ are \emph{dual}, and we write $\chi \sim \tilde{\chi}$ if for any dual bases $(u_{1},\ldots, u_{n})$ and $(\tilde{u}_{1},\ldots, \tilde{u}_{n})$, and every $1\le i\le n$, we have
\[\chi(u_{i})+\tilde{\chi}(\tilde{u}_{i})\ge 0.\]
In addition, the second-moment Lyapunov exponent $\chi$ is called
\emph{regular} if $\chi \sim \tilde{\chi}$ and $\gamma(\chi,\tilde{\chi})= 0$.

Now we illustrate that the exponents $\chi$ associated with \x{a1} and $\tilde{\chi}$
associated with \x{b2} are dual.  For this purpose, let $u(t)$ be a solution of  \x{a1}, and $\tilde{u}(t)$ be a solution of  \x{b2}. Obviously, $u(t)=\Phi(t)u_{0}$, and $\tilde{u}(t)=\Phi^{-T}(t)\tilde{u}_{0}$. Thus,
 for every $t\in I$, we have
 \beaa  \langle u(t), \tilde{u}(t)\rangle \EQ (\Phi(t)u_{0})^{T}(\Phi^{-T}(t)\tilde{u}_{0})=u_{0}^{T} \tilde{u}_{0} = \langle u_{0}, \tilde{u}_{0}\rangle,
 \eeaa
 %\beaa d \langle v(t), w(t)\rangle \EQ
%\langle A(t)v(t)dt+G(t)v(t)d\omega(t), w(t)\rangle\\
%\EM +\langle v(t),  \left[-A(t)+G^{2}(t)\right]^{*}w(t)dt-G^{*}(t)w(t)d\omega(t)\rangle\\
%\EM -\langle G(t)v(t)dt, G^{*}(t)w(t)\rangle\\
%\EQ
%\langle A(t)v(t),w(t)\rangle dt+\langle G(t)v(t), w(t)\rangle d\omega(t) -\langle A(t)v(t),w(t)\rangle dt\\
%\EM +\langle G^{2}(t)v(t),w(t)\rangle dt-\langle G(t)v(t), w(t)\rangle d\omega(t)- \langle G^{2}(t)v(t),w(t)\rangle dt \\
%\EQ 0,\eeaa
where $\langle \cdot,\cdot\rangle$ denotes the standard inner product in $\R^{n}$. Hence
\be\lb{b6}\langle u(t),\tilde{u}(t)\rangle=\langle u(0),\tilde{u}(0)\rangle\ee
for any $t\in I$. \x{a1} and \x{b2} can be called \emph{dual} due to the fact that \x{b6} holds. Now choose dual spaces $(u_{1},\ldots, u_{n})$ and $(\tilde{u}_{1},\ldots, \tilde{u}_{n})$ of $\R^{n}$. Let $u_{i}(t)$ be the unique solution of \x{a1} with $u_{i}(0)=u_{i}$, and $\tilde{u}_{i}(t)$ be the unique solution of \x{b2} with $\tilde{u}_{i}(0)=\tilde{u}_{i}$. With the help of H\"{o}lder's inequality, we have
\[\E\|u_{i}(t)\|^{2}\cdot \E\|\tilde{u}_{i}(t)\|^{2}\ge 1\]
for every $t\ge 0$, and hence, $\chi(u_{i})+\tilde{\chi}(\tilde{u}_{i})\ge 0$ for every $i$.
Thus, $\gamma(\chi,\tilde{\chi})\ge 0$  follows immediately from the analysis above.

\section{Nonuniform mean-square exponential dichotomy }

\setcounter{equation}{0} \noindent

In Section 1 we introduce the notion of NMS-ED
for SDEs,  which extends the concept of (uniform) MS-ED, and allows us to detect and formulate ``random"  versions of nonuniform behavior
  for SDEs. In this section,  we will show that {\rm \x{a1}} admits an
  NMS-ED, if there is a nonuniform Lyapunov matrix
  $S(t)$, which transforms {\rm \x{a1}} into a new system with block form.

For the convenience of later discussion, we first derive an equivalent definition of the
NMS-ED  of \x{a1}.

\begin{lemma} \lb{lemma31}
  The  projector   of   {\rm\x{a1}} can   be   chosen  as
   \[\tilde{P}=\left(\begin{array}{ll}
I_{n_{1}\times n_{1}} & 0_{n_{1}\times n_{2}}\\
0_{n_{2}\times n_{1}} & 0_{n_{2}\times n_{2}}
\end{array}\right)\]
with $n_{1}=\dim \im \tilde{P}$ and $n_{2}=\dim \ker \tilde{P}$ such that
$\tilde{P}=\Phi^{-1}(t)P(t)\Phi(t)$ hold for all $t\in I$. Thus %the fundamental matrix $\Phi(t)$ of {\rm\x{a1}} can   be
%chosen   appropriately   such   that
 the inequalities {\rm \x{c2}} can   be   rewritten   as
\be\lb{c11}
\begin{split}
\E\|\Phi(t)\tilde{P}\Phi^{-1}(s)\|^2 \leq K e^{-\alpha(t-s)+\beta s},
\quad \forall~ (t,s)\in I^{2}_{\geq},\\
\E\|\Phi(t)\tilde{Q}\Phi^{-1}(s)\|^2 \leq K e^{-\alpha(s-t)+\beta s}, \quad
\forall~ (t,s)\in I^{2}_{\leq},
\end{split}
\ee
where $\tilde{Q}=Id-\tilde{P}$.
\end{lemma}

\prf{Let $P(t)=\Phi(t)\tilde{P}\Phi^{-1}(t)$ for any $t\in I$.
Then
\[\E\|\Phi(t)\Phi^{-1}(s)P(s)\|^2=\E\|\Phi(t)\Phi^{-1}(s)\Phi(s)\tilde{P}\Phi^{-1}(s)\|^2
=\E\|\Phi(t)\tilde{P}\Phi^{-1}(s)\|^2.\]
Obviously, \x{c11} follows immediately from \x{c2}.

Conversely, it follows from \x{c1} that
\beaa P(t)\EQ  P(t)\Phi(t)\Phi^{-1}(s)\Phi(s)\Phi^{-1}(t)\\
\EQ \Phi(t)\Phi^{-1}(s)P(s)\Phi(s)\Phi^{-1}(t)
\eeaa
for any $t,~s\in I$. Then we have
\[\Phi^{-1}(t)P(t)\Phi(t)=\Phi^{-1}(s)P(s)\Phi(s)\] for all $t,~s\in I$.
Denote $\tilde{P}:=\Phi^{-1}(t)P(t)\Phi(t)$. Thus \x{c2} follows immediately from \x{c11} that
\bea\lb{fjx1} \E\|\Phi(t)\tilde{P}\Phi^{-1}(s)\|^2=\E\|\Phi(t)\Phi^{-1}(s)P(s)\Phi(s)\Phi^{-1}(s)\|^2
=\E\|\Phi(t)\Phi^{-1}(s)P(s)\|^2.\eea
Similarly to \x{fjx1}, one can prove that
\[\E\|\Phi(t)\tilde{Q}\Phi^{-1}(s)\|^2=\E\|\Phi(t)\Phi^{-1}(s)Q(s)\|^2.\]
In addition,
\[P(t)\Phi(t)\Phi^{-1}(s)=\Phi(t)\tilde{P}\Phi^{-1}(s)=\Phi(t)\Phi^{-1}(s)P(s),\]
and this completes the proof.
 \hspace{\stretch{1}}$\Box$
}

%Let $a>b>0$ be real parameters and consider the stochastic differential equation
% \be\lb{c3}dx=(-a-bt\sin t)x(t)dt+\sqrt{2b\cos t}\exp(-at+bt\cos t)d\omega(t), \quad x(0)=1. \ee
% \begin{example}The stochastic differential equation {\rm\x{c3}} admits a nonuniform exponential
%dichotomy in mean square that is not a uniform exponential dichotomy in mean square. \end{example}
%\noindent\prf{It is easy to verify that
% \[\Phi(t)=\exp\left(-at+bt\cos t-b\sin t\right)\]
%is the fundamental matrix of the equation
%\[x'=(-a-bt\sin t)x(t).\] Hence, Following the idea of \cite[pp. 97]{ev}, the solution of \x{c3} is given by
%\[x(t)=\exp\left(-at+bt\cos t-b\sin t\right)\left(1+\sqrt{2b}\int_{0}^{t}e^{b\sin s}\sqrt{\cos s}d\omega(s)\right)\]
%since $x(0)=1$. Therefore, \beaa \E\|x(t)\|^{2}\EQ \exp\left(-2at+2bt\cos t-2b\sin t\right)\left(1+2b\int_{0}^{t}e^{2b\sin s}\cos sds\right)\\
%\EQ \exp\left(-2at+2bt\cos t\right).\eeaa
%Thus \[\E\|x(t)\|^{2}=U(t,s)\E\|x(s)\|^{2}\] for every $t\ge s$, where $U(t,s)=e^{-2a(t-s)+2b(t\cos t-s\cos s)}$. Note first that
%\[U(t,s)=e^{(-2a+2b)(t-s)+2bt(\cos t-1)-2as(\cos s-1)}\] and thus \be\lb{c4}U(t,s)\le e^{(-2a+2b)(t-s)+2b s}.\ee
%Furthermore, if $t =4k\pi$ and $s =3k\pi$
%with $k\in \N$, then
%\be\lb{c5}U(t,s)= e^{(-2a+2b)(t-s)+2b s}.\ee
%Thus, the stochastic differential equation \x{c3} admits a
%nonuniform exponential dichotomy in mean square. By \x{c5}, the exponential $e^{2b s}$ in \x{c4} cannot be removed by making $a-b$ sufficiently large. This shows
%that the exponential dichotomy in mean square is not uniform.

For $x'=A(t)x$,  Barreia and Valls \cite{bv07} introduced and investigated the nonuniform behavior of $x'=A(t)x$ with the assumption that $A(t)$ has the following block form
 \[A(t)=\left(\begin{array}{ll}
A_{1}(t) & 0\\
0 & A_{2}(t)
\end{array}\right).\]
 But for \x{a1}, it is   unreasonable to assume that $A(t)$ and $G(t)$ in system \x{a1} can be decoupled into block forms with the same dimensions.
 In order to overcome the obstacle caused by the  drift term $A(t)$ and the  diffusion term $G(t)$ in \x{a1}, we construct a linear SDE  which is kinematically similar to \x{a1}. For this purpose  we  establish several auxiliary results.

Consider a linear SDE
\be\lb{c8}dv(t)=B(t)v(t)dt+v(t)d\omega(t)\ee with continuous function $B: I \rightarrow \R^{n\times n}$. Eq. \x{a1} is said to be \emph{kinematically similar} to Eq. \x{c8} if there exists a stochastic process $S(t)=(S_{ij}(t))_{n \times n}$ with
\[\sup_{t\in I}\|S(t)\|_{2}<\oo \quad {\rm and} \quad  \sup_{t\in I}\|S^{-1}(t)\|_{2}<\oo,\]
 which satisfies the stochastic differential equation
\be\lb{c9}dS(t)=(A(t)S(t)-S(t)B(t)+S(t)-G(t)S(t))dt+(G(t)S(t)-S(t))d\omega(t).\ee
The change of variables $u(t)=S(t)v(t)$ then transforms \x{a1} into \x{c8}. This technique is similar to the one used in ODEs (See e.g., \cite[p. 38]{cop} for a detailed exposition).

\begin{lemma}\lb{lemma32}
  Let  $P: I \rightarrow \R^{n\times n}$  be a  symmetric   projection, and let
$\Phi(t)$  be an invertible random matrix for any $t\in I$. The mapping
\[\tilde{R}: I \rightarrow \R^{n\times n},\quad t \rightarrow P\Phi^{T}(t)\Phi(t)P +(Id-P)\Phi^{T}(t)\Phi(t)(Id-P)\] is a positive definite, symmetric matrix for every
$t\in I$. Moreover, there exists a unique $R(t)$, $t\in I$  with
\be\lb{c24} R^{2}(t)=\tilde{R}(t), \quad PR(t)=R(t)P.
\ee
In addition, if we put
\[S: I \rightarrow \R^{n\times n},\quad t \rightarrow \Phi(t)R^{-1}(t),\]
then $S(t)$ is an invertible random matrix, which satisfies
\[S(t)PS^{-1}(t)=\Phi(t)P\Phi^{-1}(t)\]
and
\beaa \|S(t)\|^{2}_{2}\EQ \E\|S(t)\|^{2}\le 2,\\
\|S^{-1}(t)\|^{2}_{2}\EQ \E\|S^{-1}(t)\|^{2}\le \E\|\Phi(t)P\Phi^{-1}(t)\|^{2}+\E\|\Phi(t)(Id-P)\Phi^{-1}(t)\|^{2}.
\eeaa
\end{lemma}

The above lemma  is a stochastic version of estimation of kinematical similarity for ODE, which can be proved following the same way as in \cite[Lemma 1, p. 39]{cop}, so we omit the proof. One can also see Lemma A.5 in \cite{sie-02} for details.

In   the   setting   of   classical   exponential  dichotomies, $S^{-1}(t)$ is bounded,
which follows from the properties $\|\Phi(t)P\Phi^{-1}(t)\|<+\oo$ and
$\|\Phi(t)(Id-P)\Phi^{-1}(t)\|<+\oo$ (see Definition 2.1 in \cite{sie-02} for details).
 However,   in   the setting   of   NMS-ED, $S^{-1}(t)$ can be unbounded on $I$ in the nonuniform mean-square sense due to \x{c11}, i.e.,
\[\E\|\Phi(t)P\Phi^{-1}(t)\|^{2}\le K e^{\beta t}\quad {\rm and}\quad \E\|\Phi(t)(Id-P)\Phi^{-1}(t)\|^{2} \le K e^{\beta t}.\]

Based on this observation, we,  unlike the previous work in \cite{dk, sie-02}, need to consider the new notion of   nonuniform kinematical   similarity  to overcome the difficulties  arising from the lack of boundedness condition.

\begin{definition} \lb{def31}
Suppose that $S(t)$  is a stochastic process. $S(t)$ is
said   to   be  a   {\rm nonuniform  Lyapunov  matrix}   if   there  exists  a   constant $M>0$ such that \be \lb{c12} \|S(t)\|^{2}_{2}\le M e^{\beta t}, \quad {\rm and} \quad \|S^{-1}(t)\|^{2}_{2}\le M e^{\beta t}, \quad {\rm for~all ~}t\in I.
  \ee
  {\rm\x{a1}} and {\rm\x{c8}}  are said to be {\rm nonuniformly kinematically   similar} if there exists a $\R^{n\times n}$-valued invertible stochastic process  $S(t)$ satisfying  {\rm\x{c9}}.
\end{definition}

The following lemma illustrates the construction of \x{c9}. For the corresponding deterministic version of Lemma \ref{lem31}, we refer to \cite[Lemma 2.1, p. 158]{dk}.

\begin{lemma}\lb{lem31} For  a stochastic process $S(t)$, the following statements are equivalent:
\begin{enumerate}
  \item The systems {\rm\x{a1}} and {\rm\x{c8}} are nonuniformly kinematically similar via $S(t)$ on $I$;
  \item Let $\Phi_{A}(t)$ and $\Phi_{B}(t)$ denote the fundamental matrix  solutions of  {\rm\x{a1}} and {\rm\x{c8}} respectively. The identity
   \be\lb{c10}\Phi_{A}(t)\Phi^{-1}_{A}(\tau)S(\tau)=S(t)\Phi_{B}(t)\Phi^{-1}_{B}(\tau)\ee
      holds for all $t, \tau \in I$;
  \item The stochastic process  $S(t)$ solves the SDE {\rm\x{c9}}.
\end{enumerate}
 \end{lemma}

 \noindent\prf{First, assume that {\rm\x{a1}} and {\rm\x{c8}} are nonuniformly kinematically similar via $S(t)$ on $I$. Then we  obtain from $u(t)=S(t)v(t)$ the relation
 \[\Phi_{A}(t)\Phi^{-1}_{A}(\tau)u(\tau)=S(t)\Phi_{B}(t)\Phi^{-1}_{B}(\tau)v(\tau).\]
 By the arbitrariness  of $u(t)$ and the formula $u(\tau)=S(\tau)v(\tau)$, we have \[\Phi_{A}(t)\Phi^{-1}_{A}(\tau)S(\tau)=S(t)\Phi_{B}(t)\Phi^{-1}_{B}(\tau).\]
Second, assume that \x{c10} holds for all $t, \tau \in I$. Then we have
 \[\Phi_{A}(t)\Phi^{-1}_{A}(0)S(0)=S(t)\Phi_{B}(t)\Phi^{-1}_{B}(0).\]
 Denote $\tilde{\Phi}_{A}(t)=\Phi_{A}(t)\Phi^{-1}_{A}(0)S(0)$ and $\tilde{\Phi}_{B}(t)=\Phi_{B}(t)\Phi^{-1}_{B}(0)$. Hence, the operator $S(t)$ can be written as:
\[S(t)=\tilde{\Phi}_{A}(t)\tilde{\Phi}_{B}^{-1}(t).\]
It follows from \x{a1}, \x{b5} and It\^{o} product rule that
\beaa d S(t) \EQ d (\tilde{\Phi}_{A}(t)\tilde{\Phi}_{B}^{-1}(t))\\
\EQ d \tilde{\Phi}_{A}(t)\tilde{\Phi}_{B}^{-1}(t) +  \tilde{\Phi}_{A}(t)d\tilde{\Phi}_{B}^{-1}(t) +d\tilde{\Phi}_{A}(t)d\tilde{\Phi}_{B}^{-1}(t)\\
\EQ A(t)S(t)dt +G(t)S(t)d\omega(t)+S(t)(-B(t)+Id)dt-S(t)d\omega(t)-G(t)S(t)dt\\
\EQ (A(t)S(t)-S(t)B(t)+S(t)-G(t)S(t))dt+(G(t)S(t)-S(t))d\omega(t),
\eeaa
which means that Statement (3) holds true.

Finally, assuming that $S(t)$ is a fundamental matrix  solution of  SDE {\rm\x{c9}},
it follows from It\^{o} product rule that
 \beaa d(S(t)v(t))\EQ dS(t)v(t)+S(t)dv(t)+dS(t)dv(t)\\\EQ (A(t)S(t)-S(t)B(t)+S(t)-G(t)S(t))v(t)dt+(G(t)S(t)-S(t))v(t)d\omega(t)\\
 \EM + S(t)B(t)v(t)dt+S(t)v(t)d\omega(t)+(G(t)S(t)-S(t))v(t)dt \\
 \EQ A(t)S(t)v(t)dt+G(t)S(t)v(t)d\omega(t)\\
 \EQ A(t)u(t)dt+G(t)x(t)d\omega(t)=du(t).
 \eeaa
 This completes the proof of the lemma. \hspace{\stretch{1}}$\Box$

%\begin{remark}
%   Following the same idea as in dealing with the ordinary differential system (see \cite[Lemma 1, p. 39]{cop} for details).
%\end{remark}

\begin{lemma}\lb{lem32} Assuming that the systems {\rm\x{a1}} and {\rm\x{c8}} are
nonuniformly kinematically similar via $S(t)$ on $I$, and that the system {\rm\x{c8}}
admits an NMS-ED  with   the  form   {\rm\x{c11}} and
$rank(\tilde{P})=k(0\le k \le n)$, then the system {\rm\x{a1}} also admits
an NMS-ED with no change in the  projector.
 \end{lemma}

  \noindent\prf{Suppose that {\rm\x{a1}} and {\rm\x{c8}} are nonuniformly kinematically similar via $S(t)$ on $I$, and \x{c12} holds. Namely, let $\Phi_{A}(t)$ be the fundamental matrix  solution of   {\rm\x{a1}}, and $\Phi_{A}(t)=S(t)\Phi_{B}(t)$. It follows from the proof of Lemma \ref{lem31} that $\Phi_{B}(t)$ is the fundamental matrix  solution of \x{c8}. Hence, for any $t\in I$,
\bea \lb{c13} \E\|\Phi_{A}(t)\tilde{P}\Phi_{A}^{-1}(s)\|^2 \EQ \E\|S(t)\Phi_{B}(t)\tilde{P}\Phi_{B}^{-1}(s)S^{-1}(s)\|^2 \nn\\
\LE \|S(t)\|^{2}_{2} \cdot \E\|\Phi_{B}(t)\tilde{P}\Phi_{B}^{-1}(s)\|^2 \cdot \|S^{-1}(t)\|^{2}_{2} \nn\\ \LE KM^{2}e^{\beta t}e^{-\alpha(t-s)+\beta s}e^{\beta s} \nn\\
\EQ  KM^{2}e^{-(\alpha-\beta)(t-s)+3\beta s}, \quad \forall~ (t,s)\in I^{2}_{\geq}.
\eea
Similarly, one can prove that
\be \lb{c14} \E\|\Phi_{A}(t)\tilde{Q}\Phi_{A}^{-1}(s)\|^2  \le KM^{2}e^{(\alpha+\beta)(t-s)+3\beta s}, \quad \forall~ (t,s)\in I^{2}_{\leq}.
\ee
It follows from \x{c13}-\x{c14} that \x{a1} admits   an
NMS-ED due to the fact that $\beta \in [0, \alpha)$, and   there is no change in the
 projector. \hspace{\stretch{1}}$\Box$}

  \begin{lemma}\label{lem33} Assuming that  {\rm\x{a1}} admits
  an   NMS-ED of  the  form  {\rm\x{c11}} with invariant   projector $\tilde{P}\neq 0, Id$,
 the system  {\rm\x{a1}} is nonuniformly kinematically similar to a decoupled system  {\rm\x{c15}}
%\bea\lb{c15}dv(t)\EQ B(t)v(t)dt+v(t)d\omega(t)\nn\\
%\EQ\left(\begin{array}{ll}
%B_{1}(t) & 0\\
%0 & B_{2}(t)
%\end{array}\right)v(t)dt+\left(\begin{array}{ll}
%Id_{n_{1}\times n_{1}} & 0_{n_{1}\times n_{2}}\\
%0_{n_{2}\times n_{1}} & Id_{n_{2}\times n_{2}}
%\end{array}\right)v(t)d\omega(t)\eea
with
 \[B_{1}:\mathbb{R}\rightarrow \mathbb{R}^{n_{1}\times n_{1}},\quad {\rm and} \quad B_{2}:\mathbb{R}\rightarrow \mathbb{R}^{n_{2}\times n_{2}},\]
where $n_{1}=\dim \im \tilde{P}$ and $n_{2}=\dim \ker \tilde{P}$.
%possesses a projector $P: \R \rightarrow \R^{n \times n}$ such that $\|P(t)\|\le M$, $\|Id-P(t)\|\le M$ for all $t\in \R$ with a constant $M\ge 1$.

\end{lemma}
  \noindent\prf{ Let $\Phi_{A}(t)$ and $\Phi_{B}(t)$ be the fundamental matrix  solutions of {\rm\x{a1}} and {\rm\x{c8}} respectively.
  Since {\rm\x{a1}} admits an   NMS-ED of  the  form  {\rm\x{c11}} with invariant   projector $\tilde{P}\neq 0, Id$, by Lemma \ref{lemma31}, we  can   choose   a   fundamental  matrix  solution $\Phi_{A}(t)$ and the  projector $\tilde{P}=\left(\begin{array}{ll}
Id_{n_{1}\times n_{1}} & 0\\
0 & 0
\end{array}\right)(n_{1}=\dim \im \tilde{P})$ such that \x{c11} holds. For the given fundamental  matrix  solution $\Phi_{A}(t)$,
  it follows from Lemma \ref{lemma32} that there exists an invertible random matrix $S(t)=\Phi_{A}(t)\Phi^{-1}_{B}(t)$ such that
\beaa \|S(t)\|^{2}_{2}\EQ \E\|S(t)\|^{2}\le 2,\\
\|S^{-1}(t)\|^{2}_{2}\EQ \E\|S^{-1}(t)\|^{2}\le \E\|\Phi_{A}(t)\tilde{P}\Phi_{A}^{-1}(t)\|^{2}+\E\|\Phi_{A}(t)(Id-\tilde{P})\Phi_{A}^{-1}(t)\|^{2},
\eeaa
which  combined  with   the   estimates \x{c11} that
\beaa \|S(t)\|^{2}_{2}\LE 2,\\
\|S^{-1}(t)\|^{2}_{2} \LE \E\|\Phi_{A}(t)\tilde{P}\Phi_{A}^{-1}(t)\|^{2}+\E\|\Phi_{A}(t)(Id-\tilde{P})\Phi_{A}^{-1}(t)\|^{2}\le 2Ke^{\beta t}.
\eeaa
Let $M=\max\{2, 2K\}$, and we have
\[\|S(t)\|^{2}_{2}\le M e^{\beta t}, \quad {\rm and} \quad \|S^{-1}(t)\|^{2}_{2}\le M e^{\beta t}, \quad {\rm for~all ~}t\in I,\]
which implies that $S(t)$ is a nonuniform Lyapunov matrix.   Now   we  show  that   $B(t)$ has  the   block  diagonal   form  of
\x{c15}. By \x{c24}, $\Phi_{B}(t)$ commutes with matrix $\tilde{P}$ for every $t\in I$, i.e.,
\be\lb{c16} \tilde{P}\Phi_{B}(t)=\Phi_{B}(t)\tilde{P}.\ee
In addition,
\be\lb{c17}d(\Phi_{B}(t)\tilde{P})=B(t)\Phi_{B}(t)\tilde{P}dt +\Phi_{B}(t)\tilde{P}d\omega(t)
\ee
since $\Phi_{B}(t)$  is the fundamental matrix  solution of \x{c8}. By It\^{o} product rule, we have
\be\lb{c18}d(\tilde{P}\Phi_{B}(t))=\tilde{P}d\Phi_{B}(t)=\tilde{P}B(t)\Phi_{B}(t)dt +\tilde{P}\Phi_{B}(t)d\omega(t).
\ee
Taking the identity \x{c16} into \x{c17}, and comparing with \x{c18}, we have
\be\lb{c19}\tilde{P}B(t)=B(t)\tilde{P}\ee
 for every $t\in I$. Now we decompose $B: I\rightarrow \mathbb{R}^{n\times n}$ into four functions
 $$
B_{1}:I\rightarrow \mathbb{R}^{n_{1}\times n_{1}},\ B_{2}:I\rightarrow \mathbb{R}^{n_{2}\times n_{2}},
$$
$$
B_{3}:I\rightarrow \mathbb{R}^{n_{1}\times n_{2}},\ B_{4}:I\rightarrow \mathbb{R}^{n_{2}\times n_{1}},
$$
with
$$
B(t)=\left(\begin{array}{ll}
B_{1}(t) & B_{3}(t)\\
B_{4}(t) & B_{2}(t)
\end{array}\right).
$$
 Identity \x{c19} implies that

\[\left(\begin{array}{ll}
B_{1}(t) & B_{3}(t)\\
0 & 0
\end{array}\right)=\left(\begin{array}{ll}
B_{1}(t) & 0\\
B_{4}(t) & 0
\end{array}\right) \quad {\rm for}~ t\in I.\] So $B_{3}(t)\equiv 0$ and $B_{4}(t)\equiv 0$. Therefore, we get the block diagonal form
\[B(t)=\left(\begin{array}{ll}
B_{1}(t) & 0\\
0 & B_{2}(t)
\end{array}\right) \quad {\rm for}~ t\in I,\]
and the proof is complete. \hspace{\stretch{1}}$\Box$
  }

Now we can prove Theorem \ref{main1}.
%, which shows that  \x{a1} admits a nonuniform mean-square exponential dichotomy if there is a nonuniform Lyapunov matrix $S(t)$,  which transforms \x{a1} into  the block form \x{c15}

%\begin{theorem}\lb{main1}Assume that  there is a nonuniform Lyapunov matrix $S(t)$, which transforms {\rm \x{a1}} into the block form {\rm \x{c15}}.
% Then for sufficiently small $\varepsilon>0$,  {\rm \x{a1}} admits a nonuniform exponential dichotomy in mean square with
%\[\alpha =\max\{-(\chi_{k}+\varepsilon), \chi_{k+1}+\varepsilon\}>0.\]
%\end{theorem}

\noindent \mbox{\bf Proof of Theorem \ref{main1}.~}{ It suffices to prove that \x{c15}
admits an NMS-ED due to Lemma \ref{lem32} and Lemma \ref{lem33}.
From now on we consider \x{c15} with initial value $v(0)=v_{0}\in \R^{n}$.
%\be\lb{c13}dv(t)=B(t)v(t)dt+v(t)d\omega(t), \quad v(0)=v_{0} \in \R^{n}\ee.
%with $v_{0}\in \mathbb{R}^{n}$ and
%\[B(t)=\left(\begin{array}{ll}
%B_{1}(t) & 0\\
%0 & B_{2}(t)
%\end{array}\right).\]
Let $\Phi_{1}(t)$ be a fundamental matrix  solution  of the equation
\[dx(t)=B_{1}(t)x(t)dt+ x(t)d\omega(t),\] and denote by $x_{1}(t),\ldots,x_{n_{1}}(t)$
 the columns of $\Phi_{1}(t)$. Thus it follows from \x{b2} that
 $\Psi_{1}(t):=(\Phi^{-1}_{1}(t))^{T}$ is a fundamental matrix solution of the equation
\[dy(t)=\left[-B_{1}(t)+Id\right]^{T}y(t)dt-y(t)d\omega(t).\]
Also let $y_{1}(t),\ldots,y_{n_{1}}(t)$ be the columns of $\Psi_{1}(t)$. Setting
\[a_{j}=\chi(x_{j}(0))\quad {\rm and}\quad b_{j}=\tilde{\chi}(y_{j}(0))\] for each $j=1, \ldots, n_{1}$,
where $\chi$ and $\tilde{\chi}$ are the second-moment Lyapunov exponents defined as in \x{b1} and \x{b3} respectively, choosing $\varepsilon>0$ sufficiently small, there is a constant $k_{1}>1$ such that for each $j=1, \ldots, n_{1}$ and $t \in I$,
\be\lb{c20} \E\|x_{j}(t)\|^{2}\le k_{1} e^{(a_{j}+\varepsilon)t}\quad {\rm and}\quad \E\|y_{j}(t)\|^{2}\le k_{1} e^{(b_{j}+\varepsilon)t}.\ee For every $i$ and $j$, $\langle x_{i}(t), y_{j}(t)\rangle=\delta_{ij}$ follows directly from the identity $\Psi_{1}^{T}(t)\Phi_{1}(t)=Id_{n_{1}\times n_{1}}$.
In view of \x{b4}, we can assume \[\max\{a_{j}+b_{j}:j=1,\ \ldots,\ n_{1}\}=\gamma_{1},\] since the Lyapunov exponents $\chi$ and $\tilde{\chi}$ can only take a finite number of values   and the matrix $\Phi_{1}(t)$ can be chosen repeatedly until we find the minimum value. Hence the elements of the matrix $\Phi_{1}(t)\Psi_{1}^{T}(s)=\Phi_{1}(t)\Phi^{-1}_{1}(s)$ are
\[u_{ik}(t,\ s)=\sum_{j=1}^{n_{1}}x_{ij}(t)y_{kj}(s)\,\quad \forall~ (t,s)\in I^{2}_{\geq},\]
where $x_{ij}(t)$ is the $i$th coordinate of $x_{j}(t)$ , and $y_{kj}(s)$ is the $k$th coordinate of $y_{j}(s)$. Therefore,
\[
|u_{ik}(t,\ s)|^{2}\leq n_{1} \sum_{j=1}^{n_{1}}|x_{ij}(t)|^{2}\cdot|y_{kj}(s)|^{2}\leq n_{1} \sum_{j=1}^{n_{1}}\Vert x_{j}(t)\Vert^{2}\cdot\Vert y_{j}(s)\Vert^{2}.
\]
It follows from \x{c20} that
\bea\lb{c21}\E|u_{ik}(t,\ s)|^{2} \LE n_{1}\E\left(\sum_{j=1}^{n_{1}}\Vert x_{j}(t)\Vert^{2}\cdot\Vert y_{j}(s)\Vert^{2}\right) \nn\\
\LE n_{1}\sum_{j=1}^{n_{1}} \left(\E\Vert x_{j}(t)\Vert^{2} \cdot \E\Vert y_{j}(s)\Vert^{2}\right) \nn \\
\LE n_{1}k_{1}^{2}\sum_{j=1}^{n_{1}} e^{(a_{j}+\varepsilon)t+(b_{j}+\varepsilon)s}\nn\\
\EQ n_{1}k_{1}^{2}\sum_{j=1}^{n_{1}} e^{(a_{j}+\varepsilon)(t-s)+(a_{j}+b_{j}+2\varepsilon)s}\nn\\
\LE n_{1}^{2}k_{1}^{2}e^{(\chi_{k}+\varepsilon)(t-s) +(\gamma_{1}+2\varepsilon)s}, \quad \forall~ (t,s)\in I^{2}_{\geq}.
\eea
Taking $v=\displaystyle \sum_{k=1}^{n_{1}}\alpha_{k}e_{k}$ with $\displaystyle \Vert v\Vert^{2}=\sum_{k=1}^{n_{1}}\alpha_{k}^{2}=1$, where $e_{1}, \ldots, e_{n_{1}}$ is the canonical basis of $\mathbb{R}^{n_{1}}$,  we have
\bea\lb{c22}\Vert\Phi_{1}(t)\Phi^{-1}_{1}(s)v\Vert^{2}  \EQ \left\Vert\sum_{i=1}^{n_{1}}\sum_{k=1}^{n_{1}}\alpha_{k}u_{ik}(t,\ s)e_{i}\right\Vert^{2}\nn \\
\EQ \displaystyle \sum_{i=1}^{n_{1}}\left(\sum_{k=1}^{n_{1}}\alpha_{k}u_{ik}(t,\ s)\right)^{2}\leq\sum_{i=1}^{n_{1}}\left(\sum_{k=1}^{n_{1}}\alpha_{k}^{2}\sum_{k=1}^{n_{1}}u_{ik}(t,\ s)^{2}\right)\nn\\
\EQ \sum_{i=1}^{n_{1}}\sum_{k=1}^{n_{1}}u_{ik}(t,\ s)^{2}.\eea
Therefore, let $K_{1}=n_{1}^{4}k_{1}^{2}$, take  \x{c21} into \x{c22}, and we have
\bea\lb{c23}\E\Vert\Phi_{1}(t)\Phi^{-1}_{1}(s)\Vert^{2}\LE \E \left(\sum_{i=1}^{n_{1}}\sum_{k=1}^{n_{1}}u_{ik}(t,\ s)^{2}\right) \nn\\
\LE  K_{1}e^{(\chi_{k}+\varepsilon)(t-s)+(\gamma_{1}+2\varepsilon)s}, \quad \forall~ (t,s)\in I^{2}_{\geq}. \eea

Similarly, consider the matrix $\Phi_{2}(t)\Phi^{-1}_{2}(s)$, where $\Phi_{2}(t)$ is a fundamental matrix solution of the equation \[dz(t)=B_{2}(t)z(t)dt+ z(t)d\omega(t),\]
and  $\Psi_{2}(t):=(\Phi^{-1}_{2}(t))^{T}$ is a fundamental matrix  solution of the equation
\[dw(t)=\left[-B_{2}(t)+Id\right]^{T}w(t)dt-w(t)d\omega(t).\]
Let now $z_{1}(t),\ldots,z_{n_{2}}(t)$ be the columns of $\Phi_{2}(t)$, and $w_{1}(t),\ldots,w_{n_{2}}(t)$ the columns of $\Psi_{2}(t)$, and set
\[\tilde{a}_{j}=\chi(z_{j}(0))\quad {\rm and}\quad \tilde{b}_{j}=\tilde{\chi}(w_{j}(0))\] for each $j=1, \ldots, n_{2}$,
where $\chi$ and $\tilde{\chi}$ are the second-moment Lyapunov exponents defined as in \x{b1} and \x{b3} respectively. Choosing $\varepsilon>0$ sufficiently small, there is a constant $k_{2}>1$ such that for each $j=1, \ldots, n_{2}$ and $t \in I$,
\[\E\|z_{j}(t)\|^{2}\le k_{2} e^{(\tilde{a}_{j}+\varepsilon)t}\quad {\rm and}\quad \E\|w_{j}(t)\|^{2}\le k_{2} e^{(\tilde{b}_{j}+\varepsilon)t}.\] For every $i$ and $j$, $\langle z_{i}(t), w_{j}(t)\rangle=\delta_{ij}$ follows directly from the identity $\Psi_{2}^{T}(t)\Phi_{2}(t)=Id_{n_{2}\times n_{2}}$.
In view of \x{b4}, we can assume \[\max\{\tilde{a}_{j}+\tilde{b}_{j}:j=1,\ \ldots,\ n_{1}\}=\gamma_{2},\] since the Lyapunov exponents $\chi$ and $\tilde{\chi}$ can only take a finite number of values   and the matrix $\Phi_{2}(t)$ can be chosen repeatedly until we find the minimum value. Hence the elements of the matrix $\Phi_{2}(t)\Psi_{2}^{T}(s)=\Phi_{2}(t)\Phi^{-1}_{2}(s)$ are
\[v_{ik}(t,\ s)=\sum_{j=1}^{n_{2}}z_{ij}(t)w_{kj}(s)\, \quad \forall~ (t,s)\in I^{2}_{\leq},\]
where $z_{ij}(t)$ is the $i$th coordinate of $z_{j}(t)$ , and $w_{kj}(s)$ is the $k$th coordinate of $w_{j}(s)$. Therefore,
\[
|v_{ik}(t,\ s)|^{2}\leq n_{2} \sum_{j=1}^{n_{2}}|z_{ij}(t)|^{2}\cdot|w_{kj}(s)|^{2}\leq n_{2} \sum_{j=1}^{n_{2}}\Vert z_{j}(t)\Vert^{2}\cdot\Vert w_{j}(s)\Vert^{2}.
\]
Thus for all $(t,s)\in I^{2}_{\leq}$, we have
\[ \E|u_{ik}(t,\ s)|^{2}\leq n_{2}^{2}k_{2}^{2}e^{-(\chi_{k+1}+\varepsilon)(s-t) +(\gamma_{2}+2\varepsilon)s}.
\]
%Taking $v=\displaystyle \sum_{k=1}^{n_{2}}\alpha_{k}e_{k}$ with $\displaystyle \Vert v\Vert^{2}=\sum_{k=1}^{n_{2}}\alpha_{k}^{2}=1$, where $e_{1}, \ldots, e_{n_{2}}$ is the canonical basis of $\mathbb{R}^{n_{2}}$, thus we have
%\beaa\Vert\Phi_{1}(t)\Phi^{-1}_{1}(s)v\Vert^{2}  \EQ \left\Vert\sum_{i=1}^{n_{2}}\sum_{k=1}^{n_{2}}\alpha_{k}u_{ik}(t,\ s)e_{i}\right\Vert^{2}\nn \\
%\EQ \displaystyle \sum_{i=1}^{n_{2}}\left(\sum_{k=1}^{n_{2}}\alpha_{k}u_{ik}(t,\ s)\right)^{2}\leq\sum_{i=1}^{n_{2}}\left(\sum_{k=1}^{n_{2}}\alpha_{k}^{2}\sum_{k=1}^{n_{2}}u_{ik}(t,\ s)^{2}\right)\nn\\
%\EQ \sum_{i=1}^{n_{2}}\sum_{k=1}^{n_{2}}u_{ik}(t,\ s)^{2}.\eeaa
Writing $K_{2}=n_{2}^{4}k_{2}^{2}$, proceeding in a similar manner to that in \x{c22}-\x{c23},  we obtain
\beaa\E\Vert\Phi_{2}(t)\Phi^{-1}_{2}(s)\Vert^{2}\LE  K_{2}e^{-(\chi_{k+1}+\varepsilon)(s-t)+(\gamma_{2}+2\varepsilon)s}. \eeaa
 %n_{2}^{4}k_{2}^{2}e^{(\tilde{a}_{j}+\varepsilon)t+(\tilde{b}_{j}+\varepsilon)s}
%=n_{2}^{4}k_{2}^{2}e^{-(\tilde{a}_{j}+\varepsilon)(s-t)+(\tilde{a}_{j}+\tilde{b}_{j}+2\varepsilon)s}\nn\\
Therefore, we complete the proof of the theorem. \hspace{\stretch{1}}$\Box$

\section{Stability of nonuniform mean-square exponential contraction}

\setcounter{equation}{0} \noindent

We consider in this section the problems of mean-square exponential stability under the condition of NMS-EC. Eq. {\rm \x{a1}} is said to admit a \emph{nonuniform mean-square exponential contraction} if for some  constants $K, \alpha>0$ and $\beta\in [0, \alpha)$ such that  \be\lb{e14}
\E\|\Phi(t)\Phi^{-1}(s)\|^2 \leq K e^{-\alpha(t-s)+\beta |s|},
\quad \forall~ (t,s)\in I^{2}_{\geq}.
\ee
Clearly, this statement is a particular case of
NMS-ED with projection $P(t)=Id$ for every $t\in I$.
Throughout this section we assume that
$f,~ h: \R^{+}_{0}\times L^{2}(\Omega, \R^{n}) \rightarrow L^{2}(\Omega, \R^{n})$ in \x{a3} are continuous functions such that
\beaa f(t,0)=h(t,0)=0, \quad \forall t\ge 0,\eeaa
and for any $u,~v \in L^{2}(\Omega, \R^{n})$, there exist some constants $c>0$ and $q>1$ such that
\be\lb{a8}\E\|f(t,u)-f(t,v)\|^{2}\bigvee \E\|h(t,u)-h(t,v)\|^{2} \le c \E\|u-v\|^{2}(\E\|u\|^{2}+\E\|v\|^{2})^{q}\ee
for every $t\ge 0$. Here $a\vee b$ means the maximum of a and b.
The  inequality in  \x{a8} means that the perturbation in mean-square is small in the neighborhood of zero.

The following is the proof of stability result for \x{a3}.

\noindent \mbox{\bf Proof of Theorem \ref{main2}.~}{Considering the space
\beaa\mathscr{L}_{c}:=\{u: t \rightarrow
L^{2}(\Omega, \R^{n}): ~u {\rm~is ~continuous ~and~}
\|u\|_{c}\le r \}\eeaa with the norm \beaa
\|u\|_{c}=\sup\left\{(\E\|u(t)\|^2)^{\frac{1}{2}}
e^{\frac{\alpha}{2} t}:
t\ge 0\right\},\eeaa
clearly, $(\mathscr{L}_{c},
\|\cdot\|_{c})$ is a Banach spaces.
In order to state our result, we need the following lemma.
\begin{lemma}\lb{lem51}
For  any given initial value $\xi_{0} \in \R^{n}$,  the solution of Eq. {\rm \x{a3}} can be expressed as
 \bea\lb{e1} u(t)\EQ
\Phi(t)\bigg{(}\Phi^{-1}(s)\xi_{0}+ \int^{t}_{s}\Phi^{-1}(\tau)h(\tau, u(\tau))d\omega(\tau)
\nn\\
\EM + \int^{t}_{s}\Phi^{-1}(\tau)\big{(}f(\tau, u(\tau))-
G(\tau)h(\tau, u(\tau))\big{)}d\tau\bigg{)},\eea
where $\Phi(t)$ is the fundamental matrix  solution of {\rm\x{a1}} with $u(s)=\xi_{0}$.
\end{lemma}

\prf{Set
\beaa \xi(t)\EQ
\Phi^{-1}(s)\xi_{0}+ \int^{t}_{s}\Phi^{-1}(\tau)h(\tau, u(\tau))d\omega(\tau)
\\
\EM + \int^{t}_{s}\Phi^{-1}(\tau)\big{(}f(\tau, u(\tau))-
G(\tau)h(\tau, u(\tau))\big{)}d\tau.\eeaa
Clearly, $u(t)=\Phi(t)\xi(t)$,
and one can easily verify that $\xi(t)$ satisfies the differential
\beaa d \xi(t)\EQ \Phi^{-1}(t)\big{(}f(t, u(t))-
G(t)h(t, u(t))\big{)}d t +\Phi^{-1}(t)h(t, u(t))d\omega(t),~\xi(0)=\xi_{0}.
\eeaa
Since $\Phi(t)$ is a fundamental matrix solution of \x{a1}, it follows from It\^{o} product rule that
\beaa du(t)\EQ
d\Phi(t)\xi(t)+\Phi(t)d\xi(t)+G(t)\Phi(t)\Phi^{-1}(t)h(t, u(t))dt\\
\EQ A(t)u(t)dt+G(t)u(t)d\omega(t)+\big{(}f(t, u(t))-
G(t)h(t, u(t))\big{)}dt\\
\EM+h(t, u(t))d\omega(t)+G(t)h(t, u(t))dt\\\EQ
\big{(}A(t)u(t)+f(t,u(t))\big{)}dt+\big{(}G(t)u(t)+h(t,u(t))\big{)}d\omega(t),
\eeaa
which means that $u(t)=\Phi(t)\xi(t)$ is a solution of \x{a3}.  In addition, $u(s)=\xi_{0}$ is trivial, and this completes the proof of the lemma. \hspace{\stretch{1}}$\Box$
}

We proceed with the proof of Theorem \ref{main2}. In order to simplify the presentation, write $\tilde{f}(t, u(t))=f(t, u(t))-
G(t)h(t, u(t))$ in the following.
Squaring both sides of \x{e1}, and taking expectations, it follows from the
elementary inequality \be\lb{e4}\left\|\sum_{k=1}^{m}a_{k}\right\|^{2}\le
m\sum_{k=1}^{m}\|a_{k}\|^2\ee that
\bea\lb{e5} \E\|u(t)\|^{2} \LE
3\E\|\Phi(t)\Phi^{-1}(0)\xi_{0}\|^{2}+ 3\E\left\|\int^{t}_{0}\Phi(t)\Phi^{-1}(\tau)h(\tau, u(\tau))
d\omega(\tau)\right\|^{2}\nn\\
\EM + 3\E\left\|\int^{t}_{0}\Phi(t)\Phi^{-1}(\tau)\tilde{f}(\tau, u(\tau))d\tau\right\|^{2}.
\eea
We define the operator $\mathcal {T}$ in $(\mathscr{L}_{c},
\|\cdot\|_{c})$ by
\beaa (\mathcal {T}u)(t) \EQ \int^{t}_{0}\Phi(t)\Phi^{-1}(\tau)h(\tau, u(\tau))
d\omega(\tau)+\int^{t}_{0}\Phi(t)\Phi^{-1}(\tau)\tilde{f}(\tau, u(\tau))d\tau.
\eeaa
Given $u_{1}, u_{2}\in \mathscr{L}_{c}$,
it follows from \x{e5} that
\bea\lb{e8}  \E\|(\mathcal {T}u_{1})(t)-(\mathcal {T}u_{2})(t)\|^{2} \LE
3\E\left\|\int^{t}_{0}\Phi(t)\Phi^{-1}(\tau)\left[h(\tau, u_{1}(\tau))-h(\tau, u_{2}(\tau))\right]
d\omega(\tau)\right\|^{2}\nn\\
\EM + 3\E\Bigg\|\int^{t}_{0}\Phi(t)\Phi^{-1}(\tau)[\tilde{f}(\tau, u_{1}(\tau))-\tilde{f}(\tau, u_{2}(\tau))]d\tau\Bigg\|^{2}.\nn\\
\eea
On the other hand, by \x{a8}, we obtain
\bea\lb{e6}\E\|f(\tau, u_{1}(\tau))-f(\tau, u_{2}(\tau))\|^{2} \LE c \E\|u_{1}(\tau)-u_{2}(\tau)\|^{2}(\E\|u_{1}(\tau)\|^{2}+\E\|u_{2}(\tau)\|^{2})^{q}\nn\\
\LE 2^{q}c r^{2q} e^{-(q+1)\alpha \tau}\|u_{1}-u_{2}\|^{2}_{c}.
\eea
Similarly, we have
\bea\lb{e7}\E\|h(\tau, u_{1}(\tau))-h(\tau, u_{2}(\tau))\|^{2}\LE 2^{q}c r^{2q} e^{-(q+1)\alpha \tau}\|u_{1}-u_{2}\|^{2}_{c}.
\eea
By \x{e14} and \x{e7}, the first term of right-hand
side in \x{e8} can be deduced as follows:
\beaa \EM \E\left\|\int^{t}_{0}\Phi(t)\Phi^{-1}(\tau)\left[h(\tau, u_{1}(\tau))-h(\tau, u_{2}(\tau))\right]
d\omega(\tau)\right\|^{2} \\
\EQ \int^{t}_{0}\E\|\Phi(t)\Phi^{-1}(\tau)\|^{2} \E\|h(\tau, u_{1}(\tau))-h(\tau, u_{2}(\tau))\|^{2} d \tau\\
\LE 2^{q}c K r^{2q} \|u_{1}-u_{2}\|^{2}_{c}e^{-\alpha t}\int^{t}_{0}e^{-q\alpha \tau} d \tau\\
\LE \frac{2^{q}cK r^{2q} e^{-\frac{\alpha}{2} t}}{q\alpha}\|u_{1}-u_{2}\|^{2}_{c}.
\eeaa
As to the second term in \x{e8}, it follows from \x{e14}, \x{e6}, \x{e7},   $\E\|x\|\le \sqrt{\E\|x\|^{2}}$, and Cauchy-Schwarz inequality that
\beaa \EM \E\Bigg\|\int^{t}_{0}\Phi(t)\Phi^{-1}(\tau)[\tilde{f}(\tau, u_{1}(\tau))-\tilde{f}(\tau, u_{2}(\tau))]d\tau\Bigg\|^{2}\\
\EM \E\Bigg\|\int^{t}_{0}\left(\Phi(t)\Phi^{-1}(\tau)\right)^{\frac{1}{2}}
\left(\left(\Phi(t)\Phi^{-1}(\tau)\right)^{\frac{1}{2}}[\tilde{f}(\tau, u_{1}(\tau))-\tilde{f}(\tau, u_{2}(\tau))]\right)d\tau\Bigg\|^{2}\\
\LE \left(\int^{t}_{0}\E\left\|\Phi(t)\Phi^{-1}(\tau)\right\| d\tau \right)\\
\EM \times \left(\int^{t}_{0}\E\left\|\Phi(t)\Phi^{-1}(\tau)\right\| \E\left\|\tilde{f}(\tau, u_{1}(\tau))-\tilde{f}(\tau, u_{2}(\tau))\right\|^{2}d\tau \right)\\
\LE 2^{1+q}c K(1+g^{2})r^{2q} \|u_{1}-u_{2}\|^{2}_{c}\left(\int^{t}_{0}e^{-\frac{\alpha}{2}(t-\tau)} d\tau \right)\left(\int^{t}_{0}e^{-\frac{\alpha}{2}(t-\tau)}e^{-(q+1)\alpha \tau} d\tau \right)\\
\LE  \frac{2^{3+q}c K(1+g^{2})r^{2q} e^{-\frac{\alpha }{2} t}}{(2q-1)\alpha^{2}}\|u_{1}-u_{2}\|^{2}_{c}.
\eeaa
Since $q>1$, we can rewrite the inequality \x{e8} as
\beaa  \E\|(\mathcal {T}u_{1})(t)-(\mathcal {T}u_{2})(t)\|^{2}\LE \frac{2^{q+2}cK r^{2q} e^{-\frac{\alpha}{2} t}}{\alpha}\left(\frac{1}{q}+\frac{4(1+g^{2})}{(2q-1)\alpha}\right)\|u_{1}-u_{2}\|^{2}_{c}.
\eeaa
We can choose appropriate $r$  such that
\[\theta=\sqrt{\frac{2^{q+2}cK r^{2q}  }{\alpha}\left(\frac{1}{q}+\frac{4(1+g^{2})}{(2q-1)\alpha}\right)}<\frac{1}{2}.\]
Therefore,
\bea\lb{e9}  \|\mathcal {T}u_{1}-\mathcal {T}u_{2}\|_{c} \LE \theta\|u_{1}-u_{2}\|_{c}.
\eea
Given $\|\xi_{0}\|\le \delta$, and considering the operator $\tilde{\mathcal {T}}$ in $(\mathscr{L}_{c},
\|\cdot\|_{c})$ defined by
\[(\tilde{\mathcal {T}}u)(t)=\xi(t)+(\mathcal {T}u)(t)\]
with $\xi(t)=\Phi(t)\Phi^{-1}(0)\xi_{0}$, it is clear that  we have $\mathcal {T}u=0$ for $u=0$, and it follows from \x{e9} that
\[\|\mathcal {T}u\|_{c}\le \theta \|u\|_{c}.\]
On the other hand, it follows from \x{e14} that
\[\|\xi(t)\|_{c} =\sup_{t\ge 0}(\E\|\Phi(t)\Phi^{-1}(0)\xi_{0}\|^2)^{\frac{1}{2}}
e^{\frac{\alpha}{2} t}=\sqrt{K}\delta <\frac{1 }{2}r
 \]
since $\delta>0$ is sufficiently small. Therefore,
\be\lb{e10}\|\tilde{\mathcal {T}}u\|_{c}\le \|\xi(t)\|_{c}+ \|\mathcal {T}u\|_{c}\le r,\ee
and this means that  $\tilde{\mathcal {T}} \mathscr{L}_{c} \subset \mathscr{L}_{c}$. In addition, by \x{e9}, we have
\beaa \|\tilde{\mathcal {T}}u_{1}-\tilde{\mathcal {T}}u_{2}\|_{c} =\|\mathcal {T}u_{1}-\mathcal {T}u_{2}\|_{c} \LE \theta\|u_{1}-u_{2}\|_{c},
\eeaa
and thus $\tilde{\mathcal {T}}$ is a contraction in $(\mathscr{L}_{c},
\|\cdot\|_{c})$. Hence, there exists a unique  $u\in \mathscr{L}_{c}$ such that $\tilde{\mathcal {T}}u=u$.  By \x{e10}, we obtain
\[\|u\|_{c}\le \frac{1 }{2}r+\theta \|u\|_{c},\]
and thus
\[\|u\|_{c}\le \frac{r}{2(1-\theta)}.\]
Therefore the function $u(t)$ satisfies \x{e11} with $\tilde{K}=\frac{r^{2}}{4(1-\theta)^{2}}>0$. \hspace{\stretch{1}}$\Box$

\begin{remark}\lb{rem41}
  Let $\chi_{\max}$ denote the maximal value of second-moment Lyapunov exponent of {\rm\x{a1}}, and let $\gamma$ denote second-moment regularity coefficient. Using the same techniques as in the proof of  Theorem {\rm\ref{main1}}, it follows easily from {\rm\x{c23}} that $\alpha=-(\chi_{\max}+\varepsilon)$ and $\beta=\gamma+2\varepsilon$ under the condition of NMS-EC. Since $\varepsilon$ can be chosen arbitrarily small, the assumption $-q\alpha+\beta<0$ in Theorem {\rm\ref{main2}} can also be substituted by $q\chi_{\max}+\gamma<0$.
\end{remark}

From the remark above, Theorem \ref{main3} is an immediate corollary of Theorem \ref{main2}, since regularity means $\gamma=0$, and clearly, $q\chi_{\max}+\gamma<0$ implies $\chi_{\max}<0$. This is a natural  condition of NMS-EC.

\section{Second-moment Regularity Coefficient}
\setcounter{equation}{0} \noindent

Following the discussion of Remark \ref{rem41}, exponent $\alpha$ can be estimated by $\chi_{\max}$ in \x{b1} in terms of \x{a1} and its solutions.
Thus it is of special interest to derive the upper and lower bounds of the second-moment regularity coefficient $\gamma(\chi,\tilde{\chi})$, which determines the stability of  the perturbed equation. From Section 2 we know that $\gamma(\chi,\tilde{\chi}) \ge 0$.
In this section we proceed to derive the more precise lower bound and upper bound of $\gamma(\chi,\tilde{\chi})$, which have the advantage that one does not need to know any explicit information about the solutions of the linear SDE \x{a1}. More specifically, the lower and upper bounds of $\gamma(\chi,\tilde{\chi})$
can be expressed solely by the  drift term  $A(t)$, and have nothing to do with the diffusion term $G(t)$.
\subsection{Lower Bound}
\begin{theorem}
  The second-moment regularity coefficient satisfies
  \[\gamma(\chi,\tilde{\chi}) \ge \frac{2}{n}\left( \limsup_{t\rightarrow +\oo}\frac{1}{t}\int_{0}^{t} tr  A(\tau) d \tau - \liminf_{t\rightarrow +\oo}\frac{1}{t}\int_{0}^{t} tr  A(\tau)  d \tau \right).
  \]
\end{theorem}

\prf{Let $v_{1},\ldots, v_{n}$ be a basis of $\R^{n}$, and $v_{i}(t)$ be the unique solution of \x{a1} such that $v_{i}(0)=v_{i}$. Then $v_{1}(t),\ldots, v_{n}(t)$ are the columns of a fundamental matrix solution of $\Phi(t)$ of the equation \x{a1}. Thus it follows from Theorem 3.2.2 in  \cite{mao} that for every $t \ge 0$,
\be\lb{d1} \frac{\det \Phi(t)}{\det \Phi(0)}=\exp \left[ \int_{0}^{t} \left(tr  A(\tau) -\frac{1}{2} tr  G^{2}(\tau)\right) d \tau + \int_{0}^{t} tr  G(\tau)d \omega(\tau)\right].
\ee
Furthermore,
\[\E |\det \Phi(t)|^{2} \le n \prod_{j=1}^{n}\E\|v_{j}(t)\|^{2}\]
follows  directly from $|\det \Phi(t)|\le \prod_{j=1}^{n}\|v_{j}(t)\|$ and the independence of the vectors  $v_{1}(t),\ldots, v_{n}(t)$.
Thus by using  \x{d1}, log-normal distribution,
%\[ \E\left(\int_{0}^{t}g(\tau) d\omega(\tau)\right)=0, \quad \E\left(\left(\int_{0}^{t}g(\tau) d\omega(\tau)\right)^{2}\right)=\int_{0}^{t}g(\tau) d\tau, \]
we have
\beaa  \sum_{j=1}^{n}\chi(v_{j}) \EQ   \limsup_{t\rightarrow +\infty}\frac{1}{t}\log
\left( \prod_{j=1}^{n}\E\|v_{j}(t)\|^{2}\right)\\
\GE  \limsup_{t\rightarrow +\infty}\frac{2}{t} \E \left( \exp \left[ \int_{0}^{t} \left(tr  A(\tau) -\frac{1}{2} tr  G^{2}(\tau)\right) d \tau + \int_{0}^{t} tr  G(\tau)d \omega(\tau)\right]\right)\\
\EQ \limsup_{t\rightarrow +\oo}\frac{2}{t}\int_{0}^{t} tr   A(\tau)  d \tau.
\eeaa

%\[\limsup_{t\rightarrow +\oo}\frac{2}{t}\int_{0}^{t} tr   A(\tau)  d \tau \le n \sum_{j=1}^{n}\chi(v_{j}).\]

Similarly, let $w_{i}(t)$ be the unique solution of \x{b2} such that $w_{i}(0)=w_{i}$ for each $i$, where $w_{1},\ldots, w_{n}$ is another basis of $\R^{n}$. Proceeding in a similar manner, we obtain
\beaa\lb{d3}
- \sum_{j=1}^{n}\tilde{\chi}(w_{j})\LE -\limsup_{t\rightarrow +\oo}\frac{2}{t}\int_{0}^{t} \left(tr  [(-A(\tau)+G^{2}(\tau))^{T}] -\frac{1}{2} tr [(G^{T}(\tau))^{2}]-\frac{1}{2} tr [(G^{T}(\tau))^{2}]\right) d \tau \nn\\
\EQ  \liminf_{t\rightarrow +\oo}\frac{2}{t}\int_{0}^{t} \left(tr  (A(\tau)-G^{2}(\tau)) + tr (G^{2}(\tau))\right) d \tau \nn \\
\EQ \liminf_{t\rightarrow +\oo}\frac{2}{t}\int_{0}^{t} tr  A(\tau)  d \tau.
\eeaa
Therefore,
\bea \lb{d2} \EM \limsup_{t\rightarrow +\oo}\frac{2}{t}\int_{0}^{t} tr  A(\tau) d \tau - \liminf_{t\rightarrow +\oo}\frac{2}{t}\int_{0}^{t} tr  A(\tau)  d \tau \nn\\ \LE \sum_{j=1}^{n}(\chi (v_{j})+\tilde{\chi}(w_{j}) ).
\eea
Now we require that a basis $(v_{1},\ldots, v_{n})$ is dual to  a basis $(w_{1},\ldots, w_{n})$, and that the minimum in \x{b4} is obtained at this pair, i.e.,
\[\gamma(\chi,\tilde{\chi})=\max\{\chi(v_{i})+\tilde{\chi}(w_{i}):1\le i\le n\}.\]
Hence we have
\be \lb{d3} \sum_{j=1}^{n}(\chi (v_{j})+\tilde{\chi}(w_{j}) ) \le n \max\{\chi(v_{i})+\tilde{\chi}(w_{i}):1\le i\le n\} = n \gamma(\chi,\tilde{\chi}).\ee
Thus the desired result follows immediately from \x{d2} and \x{d3}.  \hspace{\stretch{1}}$\Box$
}

\subsection{Upper Bound}

\

\indent For each $k=1,\ldots,n$, denote
\[\underline{\alpha}_{k}=\liminf_{t\rightarrow +\oo}\frac{1}{t}\int_{0}^{t}a_{kk}(\tau)d\tau, \quad {\rm and}\quad \overline{\alpha}_{k}=\limsup_{t\rightarrow +\oo}\frac{1}{t}\int_{0}^{t}a_{kk}(\tau)d\tau,\]
where $a_{11}(t),\ldots,a_{nn}(t)$ are the diagonal elements of $A(t)$. In addition,
here we assume that $\overline{\alpha}_{k}, k=1,\ldots,n$ are ordered as
$\overline{\alpha}_{1}\le \ldots \le \overline{\alpha}_{n}$, since this can be
trivially achieved via row permutation of $A(t)$ and $G(t)$. Thus we
derive that the upper bound for the second-moment regularity coefficient
$\gamma(\chi,\tilde{\chi})$ can be expressed by these numbers. In the following we
use the assumption that $A(t)$ and $G(t)$
are upper triangular for every $t \in I$.
In fact,
this assumption does not affect the estimation of the upper bound of
$\gamma(\chi,\tilde{\chi})$
(see  Theorem \ref{thm53} and  Remark \ref{remark51} after the proof of Theorem \ref{thm42}).
%The reason that the assumption of $A(t)$ and $G(t)$ being upper triangular for
%every $t \in I$ does not affect the estimation of the upper bound of
%$\gamma(\chi,\tilde{\chi})$ will be explained in Theorem

%this assumption can be weakened
%by replacing the assumption of upper triangular form by the weaker assumption that \x{a1}
%can be transformed by a Lyapunov matrix into a new linear SDE with coefficient
%matrix being upper triangular form

\begin{theorem} \lb{thm42}
  Assuming that $A(t)$ and $G(t)$ are upper triangular, the second-moment regularity coefficient satisfies
\[ \gamma(\chi,\tilde{\chi})\le 2 \sum_{k=1}^{k=n}(\overline{\alpha}_{k}-\underline{\alpha}_{k}).
\]
\end{theorem}

\prf{Before proving the main result, we first present and prove several lemmas  which are useful in the proof of Theorem \ref{thm42}. The following two lemmas give  the analytic expressions of the solutions of  two kind of scalar linear SDEs.

\begin{lemma}\lb{lem41} {\rm(see \cite[Lemma 3.2.3]{mao})} Let $a(\cdot)$, $b(\cdot)$  be real-valued Borel measurable bounded functions on $[t_{0},T]$. Let
\be\lb{d5}\tilde{\Phi}(t)=e^{\int_{t_{0}}^{t}\left(a(\tau)-\frac{1}{2}b^{2}(\tau)\right)d \tau+\int_{t_{0}}^{t}b(\tau)d\omega(\tau)}.\ee
Then $x(t)=\tilde{\Phi}(t)x_{0}$ is the unique solution to the
  scalar linear SDE
\[
\left\{ \begin{array}{ll}
dx(t)  =  a(t)x(t)dt+b(t)x(t)d\omega (t),\\
x(0)  =x_{0}.
\end{array} \right.
\]
\end{lemma}

\begin{lemma}\lb{lem42} {\rm(see \cite[p. 98]{mao})} Let $a(\cdot)$, $b(\cdot)$, $c(\cdot)$ and $d(\cdot)$  be real-valued Borel measurable bounded functions on $[t_{0},T]$. $\tilde{\Phi}(t)$ is given as in {\rm\x{d5}}. Then
\[x(t)=\tilde{\Phi}(t)\left(x_{0}+\int_{t_{0}}^{t}\tilde{\Phi}^{-1}(\tau)(c(\tau)
-b(\tau)d(\tau))d\tau +\int_{t_{0}}^{t}\tilde{\Phi}^{-1}(\tau)d(\tau)d\omega (\tau)\right)\] is the unique solution to the
  scalar linear SDE
\[
\left\{ \begin{array}{ll}
dx(t)  =  (a(t)x(t)+c(t))dt+(b(t)x(t)+d(t))d\omega (t),\\
x(0)  =x_{0}.
\end{array} \right.
\]
\end{lemma}

Now we denote by $a_{ij}(t)$ and $g_{ij}(t)$ the entries of the matrix $A(t)$  and $G(t)$ respectively for each $i$ and $j$.
Denote
\[\Lambda_{i}(t):= \int_{0}^{t}\left(a_{ii}(\tau)-\frac{1}{2}g_{ii}^{2}(\tau)\right)d \tau+\int_{0}^{t}g_{ii}(\tau)d\omega(\tau),\]
and define the $n \times n$ matrix function $U(t)=(u_{ij}(t))$ as follows
\be \lb{d6} u_{ij}(t)=
\left\{ \begin{array}{lll}
0, & {\rm if}~i>j,\\
e^{\Lambda_{i}(t)},& {\rm if}~i=j,\\
\int_{0}^{t}\left(\sum_{k=i+1}^{j}a_{ik}(\tau)u_{ik}(\tau)-
g_{ij}(\tau)\sum_{k=i+1}^{j}g_{ik}(\tau)u_{ik}(\tau)\right)e^{\Lambda_{i}(t)-\Lambda_{i}(\tau)}d \tau \\
+ \int_{0}^{t}\sum_{k=i+1}^{j}g_{ik}(\tau)u_{ik}(\tau) e^{\Lambda_{i}(t)-\Lambda_{i}(\tau)}d \omega(\tau), & {\rm if}~i<j.
\end{array} \right.
\ee
 Thus it follows from Lemma \ref{lem41} and Lemma \ref{lem42} that the columns of the matrix function $U(t)$ form a basis for the space of solutions of \x{a1}.
For each $i,j=1,\ldots, n$, considering
\[\chi(u_{ij})=\limsup_{t\rightarrow +\oo}\frac{1}{t}\log \E |u_{ij}(t)|^{2},\]
we have the following result.
\begin{lemma}\lb{lem43}
  For every $i,j=1,\ldots n$, we have
\be \lb{d7} \chi(u_{ij})\le 2\left(\overline{\alpha}_{j}+\sum_{m=i}^{j-1}(\overline{\alpha}_{m}-\underline{\alpha}_{m})\right).
\ee
\end{lemma}
\prf{Firstly, it follows from \x{d6}, and log-normal distribution
%and
%\[ \E\left(\int_{0}^{t}g(\tau) d\omega(\tau)\right)=0, \quad \E\left(\left(\int_{0}^{t}g(\tau) d\omega(\tau)\right)^{2}\right)=\int_{0}^{t}g(\tau) d\tau \]
 that
\beaa\chi(u_{ii}) \EQ \limsup_{t\rightarrow +\oo}\frac{1}{t}\log \E |u_{ii}(t)|^{2}\\
\EQ \limsup_{t\rightarrow +\oo}\frac{1}{t}\log \E \left( e^{2(\int_{0}^{t}\left(a_{ii}(\tau)-\frac{1}{2}g_{ii}^{2}(\tau)\right)d \tau+2\int_{0}^{t}g_{ii}(\tau)d\omega(\tau))}\right)\\
\EQ \limsup_{t\rightarrow +\oo}\frac{2}{t}\int_{0}^{t} a_{ii}(\tau)d \tau\\
\EQ 2 \overline{\alpha}_{i}.
\eeaa

Now we can apply the backward induction method on $i$. Assuming
\be\lb{d8} \chi(u_{kj}) \le 2\overline{\alpha}_{j} +2\sum_{m=k}^{j-1}(\overline{\alpha}_{m}-\underline{\alpha}_{m}),\quad i+1\le k\le j
\ee
 holds for a given $i<n$, we  prove it for $i$.  i.e.,
\[\chi(u_{ij}) \le 2\overline{\alpha}_{j} +2\sum_{m=i}^{j-1}(\overline{\alpha}_{m}-\underline{\alpha}_{m}).\]
Clearly, for each $\varepsilon>0$, there exists $D>1$, it is easy to verify from \x{a4} and \x{d8} that
\[|a_{ik}(t)|\le D e^{t\varepsilon}, \quad \quad |g_{ik}(t)|\le D e^{t\varepsilon},\]
\[\E\left(e^{-\Lambda_{i}(t)}\right)\le D e^{(-\underline{\alpha}_{i}+\varepsilon)t},\]
and
\[\E |u_{kj}(t)|^{2}\le D e^{(2\overline{\alpha}_{j} +2\sum_{m=i+1}^{j-1}(\overline{\alpha}_{m}-\underline{\alpha}_{m})+\varepsilon)t}\]
for $t \ge 0$ and $i+1\le k\le j$. Therefore, it follows from It\^{o} isometry property, H\"{o}lder's inequality and the elementary inequality \x{e4}
that
\beaa \chi(u_{ij}) \LE \limsup_{t\rightarrow +\oo}\frac{1}{t}\log 2 \E \Biggl [ \left( \int_{0}^{t}\left(\sum_{k=i+1}^{j}a_{ik}(\tau)u_{ik}(\tau)-
g_{ij}(\tau)\sum_{k=i+1}^{j}g_{ik}(\tau)u_{ik}(\tau)\right)e^{\Lambda_{i}(t)-\Lambda_{i}(\tau)}d \tau \right)^{2} \\
\EM + \left(\int_{0}^{t}\sum_{k=i+1}^{j}g_{ik}(\tau)u_{ik}(\tau) e^{\Lambda_{i}(t)-\Lambda_{i}(\tau)}d \omega(\tau)\right)^{2}\Biggl ]\\
\LE  \limsup_{t\rightarrow +\oo}\frac{1}{t}\log \Biggl [ 4t  \int_{0}^{t}\E \left(\left(\sum_{k=i+1}^{j}a_{ik}(\tau)u_{kj}(\tau)\right)^{2}+
\left(g_{ij}(\tau)\sum_{k=i+1}^{j}g_{ik}(\tau)u_{kj}(\tau)\right)^{2}\right)\\
\EM \E \left(e^{2\Lambda_{i}(t)-2\Lambda_{i}(\tau)}\right) d \tau + 2\int_{0}^{t}\E \left(\sum_{k=i+1}^{j}g_{ik}(\tau)u_{ik}(\tau)\right)^{2} \E \left(e^{2\Lambda_{i}(t)-2\Lambda_{i}(\tau)}\right) d \tau\Biggl ]\\
\LE 2\overline{\alpha}_{i}+ \limsup_{t\rightarrow +\oo}\frac{1}{t}\log  \Biggl [ 8t \int_{0}^{t}\left(D^{4} \sum_{k=i+1}^{j} e^{(\overline{\alpha}_{j} +\sum_{m=i+1}^{j-1}(\overline{\alpha}_{m}-\underline{\alpha}_{m})
-\underline{\alpha}_{i}+4\varepsilon)\tau}\right)^2 d\tau\\
\EM + 2 \int_{0}^{t}\left(D^{3} \sum_{k=i+1}^{j} e^{(\overline{\alpha}_{j} +\sum_{m=i+1}^{j-1}(\overline{\alpha}_{m}-\underline{\alpha}_{m})
-\underline{\alpha}_{i}+3\varepsilon)\tau}\right)^2 d\tau \Biggl ]\\
\LE 2\overline{\alpha}_{i}+ \limsup_{t\rightarrow +\oo}\frac{1}{t}\log \left[\int_{0}^{t}
\left(8tD^{8}n^{2}+2D^{6}n^{2}\right)e^{(2\overline{\alpha}_{j} +2\sum_{m=i+1}^{j-1}(\overline{\alpha}_{m}-\underline{\alpha}_{m})
-2\underline{\alpha}_{i}+8\varepsilon)\tau} d \tau\right]\\
\LE  2\overline{\alpha}_{i} +2\overline{\alpha}_{j} +2\sum_{m=i+1}^{j-1}(\overline{\alpha}_{m}-\underline{\alpha}_{m})
-2\underline{\alpha}_{i}+8\varepsilon\\
\EQ 2\left(\overline{\alpha}_{j}+\sum_{m=i}^{j-1}(\overline{\alpha}_{m}-\underline{\alpha}_{m})\right)
+8\varepsilon.
\eeaa
Note that $\varepsilon>0$ is arbitrary. Thus, \x{d7} holds for every $j\ge i$, and this completes the proof of the lemma.   \hspace{\stretch{1}}$\Box$}

On the other hand, let $\tilde{A}(t):=\left(-A(t)+G^{2}(t)\right)^{T}$ and $\tilde{G}(t):=-G^{T}(t)$. Thus, it follows from \x{b2} that
\be\lb{d9} d\tilde{u}(t)=
\tilde{A}(t)\tilde{u}(t)dt+\tilde{G}(t)\tilde{u}(t)d\omega(t)\ee
is the adjoint equation of \x{a1}. Denote  the entries of the matrix $\tilde{A}(t)$  and $\tilde{G}(t)$ by $\tilde{a}_{ij}(t)$ and $\tilde{g}_{ij}(t)$ respectively for each $i$ and $j$.
Define the $n \times n$ matrix function $\tilde{U}(t)=(\tilde{u}_{ij}(t))$ as follows
\be \lb{d10} \tilde{u}_{ij}(t)=
\left\{ \begin{array}{lll}
0, & {\rm if}~i<j,\\
e^{-\Lambda_{i}(t)},& {\rm if}~i=j,\\
\int_{0}^{t}\left(\sum_{k=j}^{i-1}\tilde{a}_{ki}(\tau)\tilde{u}_{ki}(\tau)-
\tilde{g}_{ji}(\tau)\sum_{k=j}^{i-1}\tilde{g}_{ki}(\tau)\tilde{u}_{ki}(\tau)\right)e^{-\Lambda_{i}(t)+\Lambda_{i}(\tau)}d \tau \\
+ \int_{0}^{t}\sum_{k=j}^{i-1}\tilde{g}_{ki}(\tau)\tilde{u}_{ki}(\tau) e^{-\Lambda_{i}(t)+\Lambda_{i}(\tau)}d \omega(\tau), & {\rm if}~i>j.
\end{array} \right.
\ee
Thus it follows from Lemma \ref{lem41} and Lemma \ref{lem42} that the columns of the matrix function $\tilde{U}(t)$ form a basis for the space of solutions of \x{d9}.
For each $i,j=1,\ldots, n$, considering
\[\tilde{\chi}(\tilde{u}_{ij})=\limsup_{t\rightarrow +\oo}\frac{1}{t}\log \E |\tilde{u}_{ij}(t)|^{2},\]
we have the following result.
\begin{lemma} \lb{lem44}
  For every $i,j=1,\ldots n$, we have
\be \lb{d11} \tilde{\chi}(\tilde{u}_{ij})\le 2\left(-\underline{\alpha}_{j}+\sum_{m=j+1}^{i}(\overline{\alpha}_{m}-\underline{\alpha}_{m})\right).
\ee
\end{lemma}
\prf{We proceed in a similar way as the proof of Lemma  \ref{lem43}.
Firstly, it follows from \x{d6} and log-normal distribution
%and
%\[ \E\left(\int_{0}^{t}g(\tau) d\omega(\tau)\right)=0, \quad \E\left(\left(\int_{0}^{t}g(\tau) d\omega(\tau)\right)^{2}\right)=\int_{0}^{t}g(\tau) d\tau \]
 that
\beaa\tilde{\chi}(\tilde{u}_{ii}) \EQ \limsup_{t\rightarrow +\oo}\frac{1}{t}\log \E |\tilde{u}_{ii}(t)|^{2}\\
\EQ \limsup_{t\rightarrow +\oo}\frac{1}{t}\log \E \left( e^{2(\int_{0}^{t}\left(-a_{ii}(\tau)+g_{ii}^{2}(\tau)-\frac{1}{2}g_{ii}^{2}(\tau)\right)d \tau-2\int_{0}^{t}g_{ii}(\tau)d\omega(\tau))}\right)\\
\EQ \limsup_{t\rightarrow +\oo}\frac{-2}{t}\int_{0}^{t} a_{ii}(\tau)d \tau\\
\EQ 2 \underline{\alpha}_{i}.
\eeaa
Now we can apply the induction method on $i$. Assuming that
\be\lb{d12} \chi(\tilde{u}_{kj}) \le -2\underline{\alpha}_{j} +2\sum_{m=j+1}^{k}(\overline{\alpha}_{m}-\underline{\alpha}_{m}),\quad j\le k\le i-1
\ee
 holds for a given $i>1$, we  prove it for $i$.  i.e.,
\[\tilde{\chi}(\tilde{u}_{ij}) \le -2\underline{\alpha}_{j}+2\sum_{m=j+1}^{i}(\overline{\alpha}_{m}-\underline{\alpha}_{m}).\]

Clearly, for each $\varepsilon>0$, there exists $D>1$, and then one can easily verify from \x{a4} and \x{d12} that
\[|\tilde{a}_{ik}(t)|\le D e^{t\varepsilon}, \quad \quad |\tilde{g}_{ik}(t)|\le D e^{t\varepsilon},\]
\[\E\left(e^{\Lambda_{i}(t)}\right)\le D e^{(\overline{\alpha}_{i}+\varepsilon)t},\]
and
\[\E |\tilde{u}_{kj}(t)|^{2}\le D e^{(-2\underline{\alpha}_{j} +2\sum_{m=j+1}^{i-1}(\overline{\alpha}_{m}-\underline{\alpha}_{m})+\varepsilon)t}\]
for $t \ge 0$ and $j\le k\le i-1$. Therefore, it follows from It\^{o} isometry property, H\"{o}lder's inequality and the elementary inequality \x{e4}
that
\beaa \tilde{\chi}(\tilde{u}_{ij})
\LE  \limsup_{t\rightarrow +\oo}\frac{1}{t}\log \Biggl [ 4t  \int_{0}^{t}\E \left(\left(\sum_{k=j}^{i-1}\tilde{a}_{ik}(\tau)\tilde{u}_{kj}(\tau)\right)^{2}+
\left(\tilde{g}_{ij}(\tau)\sum_{k=j}^{i-1}\tilde{g}_{ik}(\tau)\tilde{u}_{kj}(\tau)\right)^{2}\right)\\
\EM \E \left(e^{-2\Lambda_{i}(t)+2\Lambda_{i}(\tau)}\right) d \tau + 2\int_{0}^{t}\E \left(\sum_{k=j}^{i-1}\tilde{g}_{ik}(\tau)\tilde{u}_{ik}(\tau)\right)^{2} \E \left(e^{-2\Lambda_{i}(t)+2\Lambda_{i}(\tau)}\right) d \tau\Biggl ]\\
\LE -2\underline{\alpha}_{i}+ \limsup_{t\rightarrow +\oo}\frac{1}{t}\log  \Biggl [ 8t \int_{0}^{t}\left(D^{4} \sum_{k=j}^{i-1} e^{(-\underline{\alpha}_{j} +\sum_{m=j+1}^{k}(\overline{\alpha}_{m}-\underline{\alpha}_{m})
+\overline{\alpha}_{i}+4\varepsilon)\tau}\right)^2 d\tau\\
\EM + 2 \int_{0}^{t}\left(D^{3} \sum_{k=j}^{i-1} e^{(-\underline{\alpha}_{j} +\sum_{m=j+1}^{k}(\overline{\alpha}_{m}-\underline{\alpha}_{m})
+\overline{\alpha}_{i}+3\varepsilon)\tau}\right)^2 d\tau \Biggl ]\\
\LE -2\underline{\alpha}_{i}+ \limsup_{t\rightarrow +\oo}\frac{1}{t}\log \left[\int_{0}^{t}
\left(8tD^{8}n^{2}+2D^{6}n^{2}\right)e^{(-2\underline{\alpha}_{j} +2\sum_{m=j+1}^{i-1}(\overline{\alpha}_{m}-\underline{\alpha}_{m})
+2\overline{\alpha}_{i}+8\varepsilon)\tau} d \tau\right]\\
\LE  -2\underline{\alpha}_{i} -2\underline{\alpha}_{j} +2\sum_{m=j+1}^{i-1}(\overline{\alpha}_{m}-\underline{\alpha}_{m})
+2\overline{\alpha}_{i}+8\varepsilon\\
\EQ 2\left(-\underline{\alpha}_{j}+\sum_{m=j+1}^{i}(\overline{\alpha}_{m}-\underline{\alpha}_{m})\right)
+8\varepsilon.
\eeaa
Note that $\varepsilon>0$ is arbitrary. Thus, \x{d11} holds for every $j\le i$, and this completes the proof of the lemma.   \hspace{\stretch{1}}$\Box$}

We now proceed with the proof of  Theorem \ref{thm42}. It follows from Lemma \ref{lem43} and  Lemma \ref{lem44} that
\[\chi(u_{j})=\max\{\chi(u_{ij}), i=1,\ldots,n\}\le 2\left(\overline{\alpha}_{j}+\sum_{m=1}^{j-1}(\overline{\alpha}_{m}-\underline{\alpha}_{m})\right),
\]
and
\[\tilde{\chi}(\tilde{u}_{j})=\max\{\chi(\tilde{u}_{ij}), i=1,\ldots,n\}\le 2\left(-\underline{\alpha}_{j}+\sum_{m=j+1}^{n}(\overline{\alpha}_{m}-\underline{\alpha}_{m})\right).
\]
Thus, we have
\be\lb{d13}\chi(u_{j})+\tilde{\chi}(\tilde{u}_{j}) \le 2 \sum_{m=1}^{n}(\overline{\alpha}_{m}-\underline{\alpha}_{m})
\ee
for every $j=1,\ldots,n$. Therefore, from the definition of the second-moment regularity coefficient $ \gamma(\chi,\tilde{\chi})$, it suffices to prove that
the bases $(u_{1},\ldots,u_{n})$ and $(\tilde{u}_{1},\ldots,\tilde{u}_{n})$ are dual. Clearly, we can let $\Phi(t)$ and $\Phi^{-T}(t)$ be fundamental matrix solutions of \x{a1} and  \x{d9} respectively. Note that the columns of the matrix function $U(t)=(u_{ij}(t))$ form the basis for the space of solutions of \x{a1}. Thus we have
 $U(t)=\Phi(t)C_{1}$ for some constant matrix $C_{1}$. Meanwhile, it is noted that the columns of the matrix function $\tilde{U}(t)=(\tilde{u}_{ij}(t))$ form the basis for the space of solutions of \x{d9}. Thus we have
 $\tilde{U}(t)=\Phi^{-T}(t)C_{2}$ for some constant matrix $C_{2}$. Therefore,
\beaa  \langle u_{i}(t), \tilde{u}_{j}(t)\rangle \EQ (U(t)u_{i}(0))^{T}(\tilde{U}(t)\tilde{u}_{j}(0))\\
\EQ (\Phi(t)C_{1}u_{i}(0))^{T}(\Phi^{-T}(t)C_{2}\tilde{u}_{j}(0))\\
\EQ(C_{1}u_{i}(0))^{T} (C_{2}\tilde{u}_{j}(0))\\
 \EQ\langle C_{1}u_{i}(0), C_{2}\tilde{u}_{j}(0)\rangle
 \eeaa
for every $t\ge 0$. In addition, it follows from \x{d6} and \x{d10} that
\beaa  u_{ij}(0)=
\left\{ \begin{array}{lll}
0, & {\rm if}~i>j,\\
1,& {\rm if}~i=j,\\
0, & {\rm if}~i<j,
\end{array} \right. \quad {\rm and} \quad \tilde{u}_{ij}(0)=
\left\{ \begin{array}{lll}
0, & {\rm if}~i<j,\\
1,& {\rm if}~i=j,\\
0, & {\rm if}~i>j.
\end{array} \right.
\eeaa
Clearly, we have $\langle u_{i}(t), \tilde{u}_{j}(t)\rangle =\langle u_{i}(0), \tilde{u}_{j}(0)\rangle=0$ for every $i\neq j$. Furthermore, it follows from \x{d6} and \x{d10} that
\beaa \langle u_{i}(t), \tilde{u}_{i}(t)\rangle = e^{\Lambda_{i}(t)}e^{-\Lambda_{i}(t)}=1.
\eeaa
Thus, it is  concluded that $\langle u_{i}(t), \tilde{u}_{i}(t)\rangle=\delta_{ij}$
for every $i$ and $j$. The theorem follows from \x{d13} and the definition of
the second-moment regularity coefficient $ \gamma(\chi,\tilde{\chi})$ immediately.
\hspace{\stretch{1}}$\Box$}

The following result implies that there exist a unitary matrix which can transform
\x{a1} into a linear SDE with coefficient
matrices being upper triangular for every $t \in I$.
\begin{theorem}\lb{thm53} There exist a unitary matrix $S(t)$ such that
the change of variable $x(t)=S^{-1}(t)u(t)$  transforms {\rm\x{a1} } into
\be\lb{d14} dx(t)=B(t)x(t)dt+H(t)x(t)d\omega(t)
\ee
with $B(t)$ and $H(t)$ being upper triangular for every $t \in I$.
\end{theorem}
\prf{Assuming that $U(t)$ is a matrix with the columns  $u_{1}(t),\ldots, u_{n}(t)$,
 where $u_{i}(t)$ is the
solution of \x{a1} satisfying the initial condition $u_{i}(0)=u_{i}$ for $i=1,\ldots,n$,
and using the Gram-Schmidt orthogonalization
procedure to the basis $u_{i}(t)$ with $i=1,\ldots,n$, we can construct
a matrix $S(t)$ with the columns  $s_{1}(t),\ldots, s_{n}(t)$ satisfying
  $\langle s_{i}(t), s_{j}(t)\rangle=\delta_{ij}$, where $\delta_{ij}$
is the Kronecker symbol. Obviously, $S(t)$  is unitary
for each $t\in I$. Moreover, the Gram-Schmidt procedure can be effected in such
a way that each function $s_{k}(t)$ is a linear combination of functions
$u_{1}(t),\ldots, u_{k}(t)$. It follows that the change of variable $X(t)=S^{-1}(t)U(t)$
is upper triangular for each $t\in I$, and the columns of $x_{1}(t)=S^{-1}(t)u_{1}(t),
\ldots, x_{n}(t)=S^{-1}(t)u_{n}(t)$ of the matrix $X(t)$ form a
basis of the space of solutions of \x{d14}.

Write $X(t)=(x_{1}(t),\ldots,x_{n}(t))$ as the following
\[X(t)=\left(
  \begin{array}{ccccc}
    x_{1,1}(t) & x_{1,2}(t) & x_{1,3}(t) & \ldots & x_{1,n}(t) \\
      & x_{2,2}(t) & x_{2,3}(t) & \ldots & x_{2,n}(t) \\
     &  & \ddots & \ddots & \vdots \\
     &   &   & \ddots & x_{n-1,n}(t) \\
    0 &   &   &   & x_{n,n}(t) \\
  \end{array}
\right),\]
since $X(t)$
is upper triangular for each $t\in I$.
  Now we prove that $B(t)$ and $H(t)$ in \x{d14} are upper triangular for each $t\in I$.
The result follows by induction. Write $B$ and $H$ in block forms:
\[ B=\left(
  \begin{array}{ccc}
    b_{11} & b_{12} & B_{13} \\
    b_{21} & b_{22} & B_{23} \\
    B_{31} & B_{32} & B_{33} \\
  \end{array}
\right), \quad {\rm and} \quad
H=\left(
  \begin{array}{ccc}
    h_{11} & h_{12} & H_{13} \\
    h_{21} & h_{22} & H_{23} \\
    H_{31} & H_{32} & H_{33} \\
  \end{array}
\right), \]
where $B_{13}, ~H_{13},~B_{23}, ~H_{23}: I \rightarrow \R^{1\times(n-2)}$,
$B_{31}, ~H_{31},~B_{32}, ~H_{32}: I \rightarrow \R^{(n-2)\times 1}$,
$B_{33}, ~H_{33}: I \rightarrow \R^{(n-2)\times (n-2)}$ are all continuous and bounded.
In order to prove that $b_{21}=h_{21}=0$ and $B_{31}=H_{31}=0^{(n-2)\times 1}$,
we choose the first column
\[x_{1}(t)=(x_{1,1}(t),\underbrace{0,\ldots,0}_{n-1})^{T}\]
of
the matrix $X(t)$,
that is
\[d \left(
               \begin{array}{c}
                 x_{11} \\
                 0 \\
                 \vdots \\
                 0\\
               \end{array}
             \right)=
\left(
  \begin{array}{ccc}
    b_{11} & b_{12} & B_{13} \\
    b_{21} & b_{22} & B_{23} \\
    B_{31} & B_{32} & B_{33} \\
  \end{array}
\right)\left(
               \begin{array}{c}
                 x_{11} \\
                 0 \\
                 \vdots \\
                 0\\
               \end{array}
             \right) d t+
\left(
  \begin{array}{ccc}
    h_{11} & h_{12} & H_{13} \\
    h_{21} & h_{22} & H_{23} \\
    H_{31} & H_{32} & H_{33} \\
  \end{array}
\right)\left(
               \begin{array}{c}
                 x_{11} \\
                 0 \\
                 \vdots \\
                 0\\
               \end{array}
             \right) d \omega(t).
\]
For the second equation of above equality, we have
\[ d 0 = b_{21}(t) x_{11}(t) dt + h_{21}(t) x_{11}(t) d \omega(t), \]
which implies that $b_{21}(t)=h_{21}(t)=0$ since $x_{1}$ is a stochastic process and $x_{11}\neq 0$.
Moreover, following the same steps as above, we obtain $B_{31}=H_{31}=0^{(n-2)\times 1}$.

Now we assume that the matrix functions $B$ and $H$ have been progressively
upper  triangulated in its first $p-1$ columns so that the transformed coefficient matrices
$B$ and $H$  have the forms
\[ B=\left(
  \begin{array}{ccc}
    B_{11} & B_{12} & B_{13} \\
    0 & b_{p,p} & B_{23} \\
    0 & B_{32} & B_{33} \\
  \end{array}
\right), \quad {\rm and} \quad
H=\left(
  \begin{array}{ccc}
    H_{11} & H_{12} & H_{13} \\
    0 & h_{p,p} & H_{23} \\
    0 & H_{32} & H_{33} \\
  \end{array}
\right), \]
where $B_{11}, ~H_{11}: I \rightarrow \R^{(p-1)\times (p-1)}$ are upper triangular, and
$B_{12}, ~H_{12}: I \rightarrow \R^{(p-1)\times 1}$,
$B_{13}, ~H_{13}: I \rightarrow \R^{(p-1)\times (n-p)}$,
$B_{23}, ~H_{23}: I \rightarrow \R^{1\times (n-p)}$,
$B_{32}, ~H_{32}: I \rightarrow \R^{(n-p)\times 1}$,
$B_{33}, ~H_{33}: I \rightarrow \R^{(n-p)\times (n-p)}$ are all continuous and bounded.
Now we prove that $B_{32}=H_{32}=0^{(n-p)\times 1}$.
To obtain this, we choose the $p$th column
\[x_{p}(t)=(x_{p,1}(t),\ldots, x_{p,p}(t), \underbrace{0,\ldots,0}_{n-p})^{T}\]
 of
the matrix $X(t)$,
that is
\[d \left(
      \begin{array}{c}
        x_{p,1} \\
        \vdots \\
        x_{p,p} \\
        0 \\
        \vdots  \\
        0 \\
      \end{array}
    \right)
=
\left(
  \begin{array}{ccc}
    B_{11} & B_{12} & B_{13} \\
    0 & b_{p,p} & B_{23} \\
    0 & B_{32} & B_{33} \\
  \end{array}
\right)\left(
      \begin{array}{c}
        x_{p,1} \\
        \vdots \\
        x_{p,p} \\
        0 \\
        \vdots  \\
        0 \\
      \end{array}
    \right) d t+
\left(
  \begin{array}{ccc}
    H_{11} & H_{12} & H_{13} \\
    0 & h_{p,p} & H_{23} \\
    0 & H_{32} & H_{33} \\
  \end{array}
\right)\left(
      \begin{array}{c}
        x_{p,1} \\
        \vdots \\
        x_{p,p} \\
        0 \\
        \vdots  \\
        0 \\
      \end{array}
    \right) d \omega(t),
\]
where we require that the entries $x_{p, j}(t)=0$, $j=p+1,\ldots, n$ satisfy
\[ d \left(
        \begin{array}{c}
          0 \\
          \vdots \\
          0 \\
        \end{array}
      \right)
= B_{32}(t) x_{p,p}(t) dt + H_{32}(t) x_{p,p}(t) d \omega(t), \]
which implies that $B_{32}=H_{32}=0^{(n-p)\times 1}$ since $x_{p}$ is a stochastic process
and $x_{p,p}\neq 0$,
the result follows. \hspace{\stretch{1}}$\Box$

}

\begin{remark}\lb{remark51}
  The assumption that $A(t)$ and $G(t)$ are upper triangular for every $t \in I$
  in Theorem {\rm\ref{thm42}} does not affect the estimation of the upper bound
  for the second-moment regularity coefficient $\gamma(\chi,\tilde{\chi})$.
 % In Theorem {\rm \ref{thm42}} one can weaken the assumption of upper triangular form
%  provided that {\rm\x{a1}}
%can be transformed by a Lyapunov matrix into a new linear SDE with coefficient
%matrix being upper triangular form.
  In fact it follows from Theorem {\rm\ref{thm53}} that  there exist a unitary matrix
  $S(t)$ such that
the change of variable $x(t)=S^{-1}(t)u(t)$  transforms {\rm\x{a1} } into {\rm \x{d14}}
 with $B(t)$ and $H(t)$ being upper triangular for every $t \in I$.  Thus
    one can follow the same idea in Lemma  {\rm\ref{lem31}} to prove that
%  the change of variable $u(t)=S(t)v(t)$ (here we assume that $\|S(t)\|=1$, if not, we can
%  make $S(t)$ by $S(t)/\|S(t)\|$) transforms {\rm\x{a1} } into
%%\be\lb{d14} dv(t)=B(t)v(t)dt+H(t)v(t)d\omega(t)
%%\ee
%with $B(t)$ and $H(t)$ being upper triangular for every $t \in I$, where
 $S(t)$ satisfies the SDE
\bea\lb{d15}dS(t)\EQ(A(t)S(t)-S(t)B(t)+S(t)H^{2}(t)-G(t)S(t)H(t))dt\nn\\
\EM+(G(t)S(t)-S(t)H(t))d\omega(t).\eea
Let $\Phi_{B,H}(t)$ be a fundamental matrix  solution of {\rm\x{d14}}, and use $\Phi_{A,G}(t):=S(t)\Phi_{B,H}(t)$ to denote the fundamental matrix  solution of {\rm\x{a1}}.
Meanwhile, one can also use the change of variable
$\tilde{x}(t)=\tilde{S}^{-1}(t)\tilde{u}(t)$ to transform
\bea \lb{d16} d\tilde{u}(t)\EQ \tilde{A}(t)\tilde{u}(t)dt+\tilde{G}(t)\tilde{u}(t)d\omega(t)\nn\\
\hh & := & \hh \left(-A(t)+G^{2}(t)\right)^{T}\tilde{u}(t)dt-G^{T}(t)\tilde{u}(t)d\omega(t)\quad
\eea
 into
\bea \lb{d17} d\tilde{x}(t)\EQ \tilde{B}(t)\tilde{x}(t)dt+\tilde{H}(t)\tilde{x}(t)d\omega(t)\nn\\
\hh & := & \hh
\left(-B(t)+H^{2}(t)\right)^{T}\tilde{x}(t)dt-H^{T}(t)\tilde{x}(t)d\omega(t).
\eea
Let  $\Phi_{\tilde{B},\tilde{H}}(t)$ be a fundamental matrix  solution of {\rm\x{d17}}, and use $\Phi_{\tilde{A},\tilde{G}}(t):=\tilde{S}(t)\Phi_{\tilde{B},\tilde{H}}(t)$ to denote the fundamental matrix  solution of {\rm\x{d16}}. It follows from Lemma {\rm\ref{lem34}} that
\[\Phi_{\tilde{A},\tilde{G}}(t)=\Phi_{A,G}^{-T}(t), \quad {\rm and}\quad
\Phi_{\tilde{B},\tilde{H}}(t)=\Phi_{B,H}^{-T}(t).\]
Thus,
\[\tilde{S}(t)=\Phi_{\tilde{A},\tilde{G}}(t)\Phi_{\tilde{B},\tilde{H}}^{-1}(t)
=\Phi_{A,G}^{-T}(t)\Phi_{B,H}^{T}(t)=S^{-T}(t).\]
Let $u(t)$ be a solution of the equation {\rm\x{a1}}, and $\tilde{u}(t)$ be
a solution of the dual equation {\rm\x{d16}}. Obviously, $x(t)=S^{-1}(t)u(t)$ and
$\tilde{x}(t)=S^{T}(t)\tilde{u}(t)$ are solutions of {\rm\x{d14}} and {\rm\x{d17}}
respectively. Hence,
 for every $t\in I$, we have
 \beaa  \langle u(t), \tilde{u}(t)\rangle \EQ (S(t)x(t))^{T}(S^{-T}(t)\tilde{x}(t))
 =x^{T}(t) \tilde{x}(t) = \langle x(t), \tilde{x}(t)\rangle,
 \eeaa
and this means that the change of variables does not affect the inner product.
%In addition, if we substitute the variable $S(t)$ with $S_{1}(t):=S(t)/\|S(t)\|$
%in {\rm\x{d15}}. We can also use the change of variable $u(t)=S_{1}(t)v(t)$
%to transform {\rm\x{a1} } into linear SDE {\rm\x{d14}} with
%the drift term $B(t)$ and the diffusion term $H(t)$
%are both upper triangular for every $t \in I$.
%Meanwhile, the variable $\tilde{S}(t)$ can also be replaced by
% $\tilde{S}_{1}(t):=\tilde{S}(t)/\|\tilde{S}(t)\|$.
Moreover, the second-moment Lyapuov exponents associated with  {\rm\x{d14}} and {\rm\x{d17}}
coincide with the second-moment Lyapuov exponents $\chi$ and $\tilde{\chi}$
associated with {\rm\x{a1}}  and {\rm\x{d16}} respectively since
$S(t)$ is unitary for each $t\in I$.  This means that the
second-moment regularity coefficient of  {\rm\x{d14}} and {\rm\x{d17}} is the same as
that of  {\rm\x{a1}}  and {\rm\x{d16}}.
Thus we can use the assumption that $A(t)$ and $G(t)$ are upper triangular for every
 $t \in I$ to compute the upper bound for the second-moment regularity coefficient
 $\gamma(\chi,\tilde{\chi})$.
\end{remark}

\section{Examples}
\setcounter{equation}{0} \noindent
%In   what   follows   we  use  an  example   to  demonstrate our result.
The following example is on the stability theory of SDE. For the perturbation of a linear SDE,
NMS-EC is not enough to guarantee the
second-moment exponential stability of its nonlinear perturbation.
This example is established by using the ideas of Perron \cite[p. 705-706]{per},
where    the  nonuniformity   arises  from   the
dependence   on  the  initial   time $s$.

%See also \cite{bp02, zlz-17} for details.

\begin{example}\lb{exam61}
  Let
  \bea\lb{f13} 0<b<a<(2 e^{-\pi}+1)b \quad {\rm and}\quad 0<\lambda<\frac{2b}{a-b}-e^{\pi}. \eea
 The following linear SDE
 \begin{equation}\lb{f8}
\left\{ \begin{array}{ll}
du_{1} & =(-a-b(\sin \log t+\cos \log t))u_{1}dt+\frac{1}{\lambda +1}u_{1}d\omega(t)\\
dv_{1} & =(-a+b(\sin \log t+\cos \log t))v_{1}dt+v_{1}d\omega(t)
\end{array} \right.
\end{equation}
admits an NMS-EC. However, any nontrivial solution of the following perturbation equation
\begin{equation}\lb{f7}
\left\{ \begin{array}{ll}
du_{2} & =(-a-b(\sin \log t+\cos \log t))u_{2}dt+\frac{1}{\lambda +1}u_{2}d\omega(t)\\
dv_{2} & =((-a+b(\sin \log t+\cos \log t))v_{2}+u_{2}^{\lambda +1})dt+v_{2}d\omega(t)
\end{array} \right.
\end{equation}
is not mean-square exponentially stable.
\end{example}
\prf{Let
\[\Phi(t)=\left(
  \begin{array}{cc}
    U(t) & 0 \\
    0 & V(t) \\
  \end{array}
\right)\]
be a fundamental matrix  solution of \x{f8}. Thus  it follows from Lemma \ref{lem41} that $u_{1}(t)=U(t)U^{-1}(s)u_{2}(s)$ and $v_{1}(t)=V(t)V^{-1}(s)v_{1}(t_{0})$  is the unique solution of \x{f8} such that
\beaa
U(t)U^{-1}(s)\EQ e^{\int_{s}^{t}-\left(a+b(\sin \log \tau+\cos \log \tau)+\frac{1}{2(\lambda +1)}\right)d\tau+\frac{1}{\lambda +1}\int_{s}^{t}d\omega(\tau)}\\
\EQ e^{-b(t\sin \log t- s\sin \log s)-\left(a+\frac{1}{2(\lambda +1)}\right)(t-s)+\frac{1}{\lambda +1}\int_{s}^{t}d\omega(\tau) },
\eeaa
and
\beaa V(t)V^{-1}(s)\EQ e^{\int_{s}^{t}(-a+b(\sin \log t+\cos \log t)-\frac{1}{2})d\tau+\int_{s}^{t}d\omega(\tau)}\\
\EQ e^{b(t\sin \log t- s\sin \log s)-(a+\frac{1}{2})(t-s)+\int_{s}^{t}d\omega(\tau)},
\eeaa
and this implies that
\bea\lb{f9} \E\|U(t)U^{-1}(s)\|^{2}\EQ e^{-2b(t\sin \log t- s\sin \log s)-2a(t-s)}\nn\\
\EQ e^{(-2a+2b)(t-s)-2bt(\sin \log t+1)+2bs(\sin\log s+1)}\nn\\
\LE e^{(-2a+2b)(t-s)+2bs},
\eea
and
\bea\lb{f10} \E\|V(t)V^{-1}(s)\|^{2}\EQ  e^{2b(t\sin \log t- s\sin \log s)-2a(t-s)}\nn\\
\EQ e^{(-2a+2b)(t-s)+2bt(\sin \log t-1)-2bs(\sin\log s-1)}\nn\\
\LE e^{(-2a+2b)(t-s)+2bs}
\eea
for all $t\ge s$.
Furthermore, if $t =e^{2k\pi+\frac{3}{2}\pi}$ and $s =e^{2k\pi+\frac{1}{2}\pi}$
with $k\in \N$, then
\be\lb{f11}\E\|U(t)U^{-1}(s)\|^{2}= e^{(-2a+2b)(t-s)+2b s}.\ee
Similarly, if $t =e^{2k\pi+\frac{1}{2}\pi}$ and $s =e^{2k\pi-\frac{1}{2}\pi}$
with $k\in \N$, then
\be\lb{f12}\E\|V(t)V^{-1}(s)\|^{2}= e^{(-2a+2b)(t-s)+2b s}.\ee
Thus,  \x{f8} admits an
NMS-EC since $-2a+2b<0$. By \x{f11} and/or \x{f12}, the exponential $e^{2b s}$ in \x{f9} and/or \x{f10} cannot be removed. This shows
that the mean-square exponential contraction  is not uniform.

%Hence, the values of the second-moment Lyapunov exponent associated with \x{f8} are
%$\chi_{1}(u_{2})=\chi_{2}(v_{2})=-2a+2b<0$.
In addition, it follows from Lemma \ref{lem51} that for any initial condition $(u_{2}(t_{0}),v_{2}(t_{0}))$,
the solution of \x{f7} is given by
\beaa
u_{2}(t)  \EQ e^{-b(t\sin \log t- t_{0}\sin \log t_{0})-\left(a+\frac{1}{2(\lambda +1)}\right)(t-t_{0})+\frac{1}{\lambda +1}\int_{t_{0}}^{t}d\omega(\tau)}u_{2}(t_{0}),\eeaa
and
\beaa
v_{2}(t)  \EQ e^{b(t\sin \log t- t_{0}\sin \log t_{0})-(a+\frac{1}{2})(t-t_{0})+\int_{t_{0}}^{t}d\omega(\tau)}\\
\EM \times \left(v_{2}(t_{0})+u^{\lambda +1}_{2}(t_{0})\int_{t_{0}}^{t} e^{-(\lambda +2)b(\tau\sin \log \tau- t_{0}\sin \log t_{0})-\lambda a(\tau-t_{0})}d \tau \right).
\eeaa
Fix $0<\delta <\frac{\pi}{4}$, and  set
\[t'_{k}=e^{2k\pi-\frac{1}{2}\pi},\quad t_{k}=e^{2k\pi-\frac{1}{2}\pi+\delta}\]
for each $k\in \N$.
Clearly, for every $\tau \in [t'_{k},t_{k}]$ we have
\[2k\pi-\frac{1}{2}\pi\le \log \tau \le 2k\pi-\frac{1}{2}\pi+\delta,\]
and
\[(2+\lambda)b\tau \cos \delta \le -(2+\lambda)b\tau \sin \log \tau.\]
This implies that
\[\int_{t'_{k}}^{t_{k}}e^{-(\lambda +2)b\tau\sin \log \tau-\lambda a \tau}d \tau \ge
\int_{t'_{k}}^{t_{k}}e^{(\lambda +2)b\tau\cos \delta-\lambda a \tau}d \tau.\]
Write $\rho=(\lambda +2)b\cos \delta-\lambda a$, thus,
\[\int_{t_{0}}^{t_{k}}e^{(\lambda +2)b\tau\cos \delta-\lambda a \tau}d \tau
= \int_{t_{0}}^{t_{k}}e^{\rho \tau}d \tau \ge \int_{t'_{k}}^{t_{k}}e^{\rho \tau}d \tau
=\frac{1-e^{-\delta}}{\rho}e^{\rho t_{k}}.\]
Let $t^{*}_{k}=e^{2k\pi+\frac{1}{2}\pi}$. Clearly, $t^{*}_{k}=e^{\pi-\delta}t_{k}>t_{k}$. Then for $k\in \N$ sufficiently large,  we obtain
\beaa e^{b t^{*}_{k}\sin \log t^{*}_{k}}\int_{t_{0}}^{t^{*}_{k}} e^{-(\lambda +2)b\tau\sin \log \tau -\lambda a\tau}\GE e^{b t^{*}_{k}}\int_{t'_{k}}^{t_{k}}e^{\rho \tau}d \tau\\
\EQ \frac{1-e^{-\delta}}{\rho}e^{b t^{*}_{k}+\rho t_{k}}=
\frac{1-e^{-\delta}}{\rho}e^{(b +\rho e^{\delta-\pi}) t^{*}_{k}}.
\eeaa
On the other hand, we have
\beaa
\E\|v_{2}(t)\|^{2}  \EQ e^{-2a(t-t_{0})+2b(t\sin \log t- t_{0}\sin \log t_{0})}\\
\EM \times \left(v_{2}(t_{0})+u^{\lambda +1}_{2}(t_{0})\int_{t_{0}}^{t} e^{-(\lambda +2)b(\tau\sin \log \tau- t_{0}\sin \log t_{0})-\lambda a(\tau-t_{0})}d \tau \right)^{2}.
\eeaa
Thus it follows from \x{f13}
that  the second-moment Lyapunov exponent of any solution of \x{f7} satisfies
\[\chi(v_{2})\ge -2a+2b+2\rho e^{\delta-\pi}=-2a+2b+2[(\lambda +2)b\cos \delta-\lambda a] e^{\delta-\pi}>0\] if $u_{2}(t_{0})\neq 0$.
Therefore, the solution $v_{2}(t)$ is not mean-square exponentially stable. This completes the construction of the example. \hspace{\stretch{1}}$\Box$
}

\end{document}